\pdfoutput=1

\RequirePackage{silence}
\documentclass[a4paper,svgnames]{amsart}

\usepackage[hmarginratio={1:1},vmarginratio={1:1},lmargin=60.0pt,tmargin=60.0pt]{geometry}

\usepackage{float}
\usepackage{multirow,multicol}

\usepackage{booktabs}

\allowdisplaybreaks

\usepackage{subcaption}

\usepackage[T1]{fontenc}

\usepackage{latexsym,exscale,amsfonts,amssymb,mathtools}
\usepackage{amsmath,amsthm,amscd,textcomp,bbm}
\usepackage{mathrsfs,stackrel}

\usepackage{etoolbox}

\usepackage[all]{xy}

\usepackage{tikz}

\usetikzlibrary{cd}
\usetikzlibrary{decorations}
\usetikzlibrary{decorations.markings}
\usetikzlibrary{decorations.pathreplacing}
\usetikzlibrary{decorations.pathmorphing}
\usetikzlibrary{arrows.meta,shapes,positioning,matrix,calc}
\usetikzlibrary{shapes.callouts}

\tikzset{
anchorbase/.style={baseline={([yshift=-0.5ex]current bounding box.center)}},
tinynodes/.style={font=\tiny, text height=0.25ex, text depth=0.05ex},
smallnodes/.style={font=\scriptsize, text height=0.75ex, text depth=0.15ex},
crossline/.style={preaction={draw=white,line width=5.0pt,-},preaction={draw=black,line width=0.9pt,-}},
usual/.style={line width=1.0,color=black},
dot/.style={
decoration={markings,
post length=0.25mm,
pre length=0.25mm,
mark=at position #1 with {\node[circle,radius=0.2cm,inner sep=-1.5pt,color=black,fill=black]{};}
},
postaction={decorate}
},
dot/.default=1,
}

\usepackage{aliascnt}

\usepackage{ytableau}

\usepackage{xparse}

\usepackage{cite}

\usepackage{dynkin-diagrams}

\usepackage{nicematrix}

\robustify{\cite}

\usepackage{tikz-3dplot}

\DeclareMathAlphabet{\mathscrbf}{OMS}{mdugm}{b}{n}

\DeclareMathOperator{\rad}{rad}

\usepackage{array}

\newcolumntype{C}{>{$}c<{$}}

\newcolumntype{P}[1]{>{\centering\arraybackslash}p{#1}} 

\DeclarePairedDelimiter\floor{\lfloor}{\rfloor}

\usepackage{xcolor,colortbl}

\definecolor{mygray}{gray}{0.6}
\definecolor{mygraydark}{gray}{0.4}
\definecolor{mygraylight}{gray}{0.85}
\definecolor{spinach}{RGB}{46,139,87}
\definecolor{tomato}{RGB}{255,99,71}
\definecolor{orchid}{RGB}{143,40,194}
\definecolor{neon}{RGB}{77,77,255}
\definecolor{lightneon}{RGB}{110,110,255}
\definecolor{pumpkin}{RGB}{224,180,80}
\definecolor{citron}{RGB}{190,180,90}
\definecolor{lava}{RGB}{207,16,32}
\definecolor{cream}{RGB}{255,253,208}
\definecolor{verdigris}{RGB}{67,179,174}
\definecolor{Black}{RGB}{0,0,0}
\definecolor{mydarkblue}{RGB}{10,10,170}
\definecolor{darkspinach}{RGB}{20,70,20}
\definecolor{darktomato}{RGB}{155,40,30}
\definecolor{darkorchid}{RGB}{50,10,100}
\definecolor{darklava}{RGB}{150,8,16}
\definecolor{zero}{RGB}{0,0,0}
\definecolor{one}{RGB}{255,0,0}
\definecolor{two}{RGB}{0,255,0}
\definecolor{three}{RGB}{0,0,255}

\newcommand{\Xop}{%
\tikzvcenter{%
\draw[usual] (0,0) -- (0.3,0.3);
\draw[usual] (0.3,0) -- (0,0.3);
}%
}

\newcommand{\Dop}{%
\tikzvcenter{%
\draw[usual] (0,0) -- (0.4,0);
\draw[usual] (0.4,0) -- (0.4,0.4);
\draw[usual] (0.4,0.4) -- (0,0.4);
\draw[usual] (0,0.4) -- (0,0);
}%
}

\newcommand{\ebot}{%
\tikzvcenter{%
\draw[usual,dot] (0,1) -- (0,0.8);
}%
}

\newcommand{\etop}{%
\tikzvcenter{%
\draw[usual,dot] (0,0) -- (0,0.2);
}%
}

\newcommand{\Cupop}{%
\tikzvcenter{%
\draw[usual] (0,0.6)
to[out=270,in=180] (0.25,0.35)
to[out=0,in=270] (0.5,0.6);
}%
}

\newcommand{\Capop}{%
\tikzvcenter{%
\draw[usual] (0,0)
to[out=90,in=180] (0.25,0.25)
to[out=0,in=90] (0.5,0);
}%
}

\usepackage{enumitem}

\setlist[enumerate]{itemsep=0.15cm,label=\emph{\upshape(\alph*)}}
\setlist[enumerate,2]{itemsep=0.15cm,label=\emph{\upshape(\roman*)}}
\setlist[enumerate,3]{itemsep=0.15cm,label=\emph{\upshape(\Alph*)}}

\newcommand{\N}{\mathbb{N}}
\newcommand{\Z}{\mathbb{Z}}

\newcommand{\R}{\mathbb{R}}
\newcommand{\C}{\mathbb{C}}

\newcommand{\tikzvcenter}[1]{%
\vcenter{\hbox{\tikz[baseline=(current bounding box.center)]{#1}}}%
}

\DeclareMathOperator{\Ind}{Ind}
\DeclareMathOperator{\Res}{Res}
\DeclareMathOperator{\Hom}{Hom}

\def\NewTheorem#1{%
\newaliascnt{#1}{equation}%
\newtheorem{#1}[#1]{#1}%
\aliascntresetthe{#1}%
\expandafter\def\csname #1autorefname\endcsname{#1}%
}

\def\equationautorefname~#1\null{(#1)\null}

\numberwithin{equation}{subsection}

\NewTheorem{Proposition}
\NewTheorem{Theorem}
\NewTheorem{Corollary}
\NewTheorem{Lemma}
\NewTheorem{Conjecture}
\NewTheorem{Speculation}
\NewTheorem{Observation}
\NewTheorem{Assumption}

\theoremstyle{definition}

\NewTheorem{Definition}
\AtEndEnvironment{Definition}{\null\hfill$\Diamond$}%

\NewTheorem{Classification Problem}
\AtEndEnvironment{Classification Problem}{\null\hfill$\Diamond$}%

\NewTheorem{Notation}
\AtEndEnvironment{Notation}{\null\hfill$\Diamond$}%

\NewTheorem{Example}
\AtEndEnvironment{Example}{\null\hfill$\Diamond$}%

\NewTheorem{Examples}
\AtEndEnvironment{Examples}{\vskip-10mm\null\hfill$\Diamond$}%

\theoremstyle{remark}

\NewTheorem{Remark}
\AtEndEnvironment{Remark}{\null\hfill$\Diamond$}%

\NewTheorem{Question}
\AtEndEnvironment{Question}{\null\hfill$\Diamond$}%

\usepackage[hypertexnames=false]{hyperref}

\usepackage{bookmark}

\hypersetup{
pdftoolbar=true,
pdfmenubar=true,
pdffitwindow=false,
pdfstartview={FitH},
pdftitle={Affine diagram categories, algebras and monoids},
pdfauthor={David He and Daniel Tubbenhauer},
pdfsubject={},
pdfcreator={David He and Daniel Tubbenhauer},
pdfproducer={David He and Daniel Tubbenhauer},
pdfkeywords={},
pdfnewwindow=true,
colorlinks=true,
linkcolor=mydarkblue,
citecolor=teal,
filecolor=magenta,
urlcolor=orchid,
linkbordercolor=lava,
citebordercolor=teal,
urlbordercolor=orchid,
linktocpage=true
}

\def\makeautorefname#1#2{\csdef{#1autorefname}{#2}}

\makeautorefname{section}{Section}%
\makeautorefname{subsection}{Section}%
\makeautorefname{subsubsection}{Section}%
\begin{document}

\arrayrulewidth=0.5mm
\setlength{\arrayrulewidth}{0.5mm}

\title[Affine diagram categories, algebras and monoids]{Affine diagram categories, algebras and monoids}
\author[D. He and D. Tubbenhauer]{David He and Daniel Tubbenhauer}
\address{D.H.: Department of Mathematics, University of Oregon, Fenton Hall, Eugene, OR 97403, USA}
\email{davidhe@uoregon.edu}

\address{D.T.: The University of Sydney, School of Mathematics and Statistics F07, Office Carslaw 827, NSW 2006, Australia, \href{http://www.dtubbenhauer.com}{www.dtubbenhauer.com}, \href{https://orcid.org/0000-0001-7265-5047}{ORCID 0000-0001-7265-5047}}
\email{daniel.tubbenhauer@sydney.edu.au}

\begin{abstract}
We introduce and study several affine (=annular in this paper) versions of 
the classical diagram algebras such as Temperley--Lieb, partition, Brauer, Motzkin, rook Brauer, rook, planar partition, and planar rook algebras. We give generators and relation{s} presentations for them and their associated categories, study their representation theory, and the asymptotic behavior of tensor products of their representations in the monoid case. Under a mild hypothesis, we also prove a previous conjecture concerning the asymptotic growth of the number of indecomposable summands in the tensor powers of representations for finite monoids.
\end{abstract}


\addtocontents{toc}{\protect\setcounter{tocdepth}{1}}

\maketitle

\tableofcontents

\section{Introduction}

In this paper we study the representation theory of affine diagram algebras, that is, suitable subalgebras of an affinization of the partition algebra.

\subsection{Diagram algebras}

Diagram algebras such as the Temperley--Lieb, Brauer, and partition algebras arise naturally in Schur--Weyl dualities as centralizer algebras of (quantum) group actions on tensor powers. The classical example is the work of Rumer--Teller--Weyl, who gave a diagrammatic description of the representation theory of $\mathrm{SL}_2$ and thereby invented what is now called the Temperley--Lieb calculus \cite{RuTeWe-sl2}. Their description shows that the representation theory of the Temperley--Lieb algebra controls tensor product multiplicities for $\mathrm{SL}_2$-representations. In other words, the internal structure of these diagram algebras reflects, in a surprisingly efficient way, the structure of the underlying tensor category. This already suggests that the representation theory of diagram algebras is both important and intrinsically rich.

The representation theory of diagram algebras has therefore been studied extensively; see for instance \cite{halverson2020set,Tubbenhauer_2024,Scrimshaw_2023} for overviews and entry points into the literature. Beyond their role in Schur--Weyl duality, these algebras appear in a wide range of contexts: statistical mechanics \cite{1971TL,PASQUIER1990523, Levy1991, jones1994potts, martin2008diagram,garbali2017dilutetemperleyliebon1loop}; knot theory \cite{jones1985,abramsky2007temperley,Khovanov_2022}; monoidal category theory \cite{hu2019presentationsdiagramcategories, mousaaid2021affinization,east2023presentations}; modular representation theory via tilting modules \cite{AnStTu-cellular-tilting,SuTuWeZh-mixed-tilting}; logic \cite{dosen2012syntaxsplitpreorders}; semigroup theory \cite{auinger2012krohn, DOLINKA2015119, Dolinka_2017,EAST201763,durdev2019,east2021ehresmanntheorypartitionmonoids}; and even cryptography \cite{KhSiTu-monoidal-cryptography,stewart2025representation,arms2025representation,liu}. These references are meant as representative samples only: there are far too many papers on diagram algebras to list them all. In fact, the algebras, their diagrammatics, and their representation theory have been rediscovered independently many times, which is perhaps the strongest indication of how ubiquitous they have become.

\subsection{Affine diagram algebras}

Given an algebra defined by diagrams, it is then natural to ask what happens when the diagrams are placed on a cylinder (or annulus); this procedure is often called affinization.

The idea of affinization of diagram algebras has been around for decades. Early occurrences include Brieskorn's description of affine braid groups \cite{MR422674,MR1911508}, skein theory on a torus (instead of an annulus) \cite{MR964255}, and annular subfactors \cite{MR1309131,MR1929335}. The affinization of monoidal categories is also well known \cite{MR1929335}, and \cite{mousaaid2021affinization} provides a useful modern reference that we will rely on below.

For example, the algebras obtained by replacing ordinary Temperley--Lieb diagrams with their cylindrical counterparts are known (among other names) as affine or periodic Temperley--Lieb algebras. They play an important role in statistical mechanics \cite{PASQUIER1990523,Levy1991,Jacobsen2009, Pinet_2022}, and their representation theory has been studied in depth, see e.g. \cite{graham1998representation, green2023representations, erdmann1998representationsaffinetemperleyliebalgebras, langlois2023uncoiled}. Just as the Temperley--Lieb algebras are quotients of Hecke algebras of type $A$, the periodic Temperley--Lieb algebras arise as quotients of Hecke algebras of type $\tilde{A}$ \cite{graham1998representation}. A version of the affine Brauer algebra appears in \cite{goodman2004affine}. The corresponding Schur--Weyl duals often involve infinite-dimensional representation theory, see for instance \cite{IoLeZh-verma-schur-weyl}. Taken together, these examples show that affine (or annular) diagram algebras are just as pervasive as their classical counterparts, but their structure and representation theory are still far from being systematically understood.

\begin{Remark}
The terminology for ``affine'' in the diagrammatic literature is far from uniform. Essentially the same objects appear under the names affine, annular, periodic, or cylindric Temperley--Lieb, Brauer, partition, etc.\ algebras, sometimes with presentational differences. We will uniformly speak of affine (or annular) diagram algebras in this paper.
\end{Remark}

\subsection{What we do}

The paper is organized as follows. In \autoref{definition section} we introduce several versions of infinite-dimensional affine diagram algebras, obtained as affinizations of the ordinary diagram algebras in the sense that they have bases consisting of Temperley--Lieb, Motzkin, Brauer, etc.\ diagrams placed on a cylinder (or annulus). We also define categories, some of which carry a natural monoidal structure, from which the affine diagram algebras may be recovered as endomorphism algebras. We give presentations for these algebras and for the associated monoidal categories.

In \autoref{rep section} we study the representation theory of these algebras. For those which admit an involutive sandwich cellular structure (in the sense of \cite{Tubbenhauer_2024}), we classify the simple modules and compute the dimensions of the cell modules by determining the top sets in the sandwich cell data. For the algebras which are not sandwich cellular, we pass to suitable finite-dimensional quotients which are sandwich cellular, classify their simple modules, and then inflate these to obtain simple modules of the original infinite-dimensional algebras.

In \autoref{growth section} we study the {\textit{growth problem} (see below)} for the {(underlying monoids of)} affine diagram algebras {in the case where they are monoid algebras and so in particular are bialgebras equipped with a comultiplication, allowing us to define tensor products of modules}. If $M$ is a monoid and $k$ a field and $V$ is a finite-dimensional $kM$-module, we associate to $V$ two numerical invariants
\[
l(n)=l^{M,V}(n)=\#\{\text{composition factors of } V^{\otimes n}\text{ (with multiplicity)}\},\] 
\[b(n)=b^{M,V}(n)=\#\{\text{$M$-indecomposable summands of } V^{\otimes n}\text{ (with multiplicity)}\}.
\]
The study of the asymptotic behavior of these statistics as $n\to\infty$ is called the growth problem associated to $V$ and has been recently studied in \cite{he2025growthproblemsrepresentationsfinite2,he2025tensorpowersrepresentationsdiagram}. For {(the underlying monoids of) the} affine diagram algebras we introduce, we determine the growth rates of $l(n)$ and $b(n)$ for a range of natural modules, thereby producing new examples in this general framework. Along the way, we prove (a more general version of) \cite[Conjecture 1]{he2025growthproblemsrepresentationsfinite2}, originally formulated for finite monoids, under some mild assumptions. In particular, the conjecture is shown to be true for finite regular monoids and some affine diagram monoids.

More generally, growth problems may be formulated for objects in any additive Krull--Schmidt monoidal category (see \cite{lacabanne2024asymptoticsinfinitemonoidalcategories}). The growth problem for various other categories has been studied in e.g.\ \cite{coulembier2023asymptotic, Coulembier_2023, lacabanne2023asymptotics, lacabanne2024asymptoticsinfinitemonoidalcategories, coulembier2024fractalbehaviortensorpowers, he2025growthproblemsrepresentationsfinite,larsen2024boundsmathrmsl2indecomposablestensorpowers, gruber2025growthproblemsdiagramcategories, osullivan2025growthproblemsquantumgroups}. Our results fit into this circle of questions by providing new growth phenomena arising from affine diagram algebras.
\smallskip

\noindent\textbf{Acknowledgments.}
DT acknowledges support from the ARC Future Fellowship FT230100489 and attests that all participants survived, albeit with slightly reduced sanity levels. {We would like to thank the referee for giving the manuscript a very careful reading, and for offering many helpful comments which improved the quality of this paper.}

\section{Affine diagram algebras}\label{definition section}

We now introduce affine diagram algebras and monoids, and their associated categories.

\subsection{Definitions}

\subsubsection{Algebras and monoids}\label{alg and monoid section}
We assume some familiarity with the classical diagram monoids.

\begin{Definition}
An \textbf{{affine} partition $m$-diagram} consists of an infinite horizontal strip $\R\times [0,1]$ with a finite number of horizontal lines (called non-contractible loops), with one vertex at each point in $\Z\times \{0,1\}$ and at most one edge between any two vertices, such that the diagram is invariant under left or right translation by $m$. Such a diagram determines a set partition of {the set of vertices $\Z\times \{0,1\}$}, and we say two vertices lie in the same component if they are in the same block of the corresponding set partition. We consider two diagrams to be the same if they have the same number of non-contractible loops and represent the same set partition. 
\end{Definition}

Since an affine partition $m$-diagram is invariant under translation, it is completely determined by what it looks like in the \textit{fundamental rectangle} $[0,m]\times [0,1]$. We draw an affine partition $5$-diagram in \autoref{fig: affine partition diagram} below, with the fundamental rectangle shaded. The diagram has (up to translation) five components, of size $4, 3,1,1,1$ respectively, and it has two non-contractible loops.
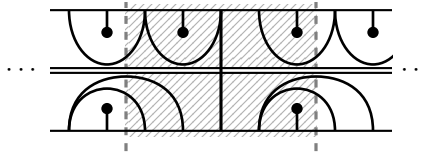
\begin{figure}[H]
\begin{center}
\begin{tikzpicture}[anchorbase]
\def\period{2.5} 
\begin{scope}
\clip (-2.25,-0.3) rectangle (2.25,1.7);
\begin{scope}[yscale=1.6]
\fill[pattern=north east lines,pattern color=black!30]
(-1.25,-0.05) rectangle (1.25,1.05);
\draw[very thick,dashed,draw=gray] (-1.25,-0.3) -- (-1.25,1.7);
\draw[very thick,dashed,draw=gray] ( 1.25,-0.3) -- ( 1.25,1.7);
\draw[thick] (-2.25,0) -- (2.25,0);
\draw[thick] (-2.25,1) -- (2.25,1);
\draw[thick, double, double distance=1pt]
(-2.25,0.5) -- (2.25,0.5);

\foreach \k in {-1,0,1} {
\begin{scope}[xshift={\k*\period cm}]
\draw[usual] (0.5,0) to[out=90,in=180] (1.25,0.45) to[out=0,in=90] (2,0);
\draw[usual] (0.5,0) to[out=90,in=180] (1,0.35)  to[out=0,in=90] (1.5,0);
\draw[usual] (0.5,1) to[out=270,in=180] (1,0.55)  to[out=0,in=270] (1.5,1);
\draw[usual] (1.5,1) to[out=270,in=180] (2,0.55)  to[out=0,in=270] (2.5,1);
\draw[usual] (0,0) -- (0,1);
\draw[usual] (2.5,0) -- (2.5,1);
\draw[usual,dot] (1,0) -- (1,0.2);
\draw[usual,dot] (1,1) -- (1,0.8);
\draw[usual,dot] (2,1) -- (2,0.8);
\end{scope}
}
\end{scope} 
\end{scope} 
\node at (-2.6,0.8) {$\cdots$};
\node at ( 2.6,0.8) {$\cdots$};
\end{tikzpicture}
\end{center}
\caption{An affine partition $5$-diagram with the fundamental rectangle shaded.}
\label{fig: affine partition diagram}
\end{figure}

Alternatively, the reader may think of these as being on a cylinder or annulus. For example, using $b=$ bottom and $t=$ top, illustrating only the fundamental rectangle:
\begin{gather*}
\begin{tikzpicture}[anchorbase,scale=1]
\fill[pattern=north east lines,pattern color=black!30]
(-0.5,-0.05) rectangle (0.5,1.05);
\draw[usual] (0,0)node[below]{b} to (-0.5,0.5);
\draw[usual] (0.5,0.5) to (0,1)node[above]{t};
\draw[usual] (-0.5,0) to (0.5,0);
\draw[usual] (-0.5,1) to (0.5,1);
\draw[very thick,dashed,draw=gray] (-0.5,-0.05) to (-0.5,1.05);
\draw[very thick,dashed,draw=gray] (0.5,-0.05) to (0.5,1.05);
\end{tikzpicture}
\leftrightsquigarrow
\begin{tikzpicture}[anchorbase,scale=1]
\draw[thick,fill=gray!50] (0,0) circle (1.5cm);
\draw[thick,fill=white] (0,0) circle (0.5cm);
\draw[usual] (0,0.5)node[below]{b} to[out=90,in=0] (-0.325,0.75) to[out=180,in=90] (-0.75,0) to[out=270,in=180] (0,-1) to[out=0,in=270] (0.75,0)to[out=90,in=270] (0,1.5)node[above]{t};
\draw[very thick,dashed,draw=gray] (0,-0.5) to (0,-1.5);
\end{tikzpicture}
\,.
\end{gather*}
Just as for ordinary partition diagrams, two affine partition $m$-diagrams could be multiplied by stacking the first diagram on top of the second one and identifying the vertices in the middle row. This produces a new affine partition $m$-diagram after removing all (contractible) closed components which arise from the stacking, and juxtaposing the non-contractible loops.   

\begin{Definition}\label{affine diagram def}
Let $k$ be a commutative ring, with $\beta \in k$. The \textbf{affine partition algebra} $\overline{aPa_m}(\beta)=\overline{aPa^k_m}(\beta)$
is the (infinite-dimensional) $k$-algebra with a basis given by  (equivalence classes) of affine partition $m$-diagrams, with multiplication defined on the basis by $D\cdot D'=\beta^{m(D,D')} D''$ and extended linearly, where $D''$ is the diagram obtained by stacking $D$ on top of $D'$
\begin{gather*}
D''=D\circ D'
=
\begin{tikzpicture}[anchorbase,smallnodes,rounded corners]
\node[rectangle,draw,minimum width=0.6cm,minimum height=0.6cm,ultra thick] at(0,0){\raisebox{-0.05cm}{$D'$}};
\node[rectangle,draw,minimum width=0.6cm,minimum height=0.6cm,ultra thick] at(0,0.6){\raisebox{-0.05cm}{$D$}};
\end{tikzpicture}
,
\end{gather*}
and $m(D,D')$ is the number of closed components removed.

We then define the following subalgebras:
\begin{enumerate}
\item The \textit{affine Brauer algebra} $\overline{aBr_m}(\beta)$ is the subalgebra of $\overline{aPa_m}(\beta)$ generated by diagrams all of whose components have size 2.
\item The \textit{affine rook-Brauer algebra} $\overline{aRoBr_m}(\beta)$ is the subalgebra of $\overline{aPa_m}(\beta)$ generated by diagrams all of whose components have size $\le 2$.
\item The \textit{affine rook algebra} $\overline{aRo_m(\beta)}=kaRo_m$ is the subalgebra of $\overline{aRoBr_m}(\beta)$ generated by diagrams all of whose components have {at most} one top vertex and one bottom vertex, and no non-contractible loops.
\item The \textit{affine Temperley--Lieb algebra} $\overline{aTL_m(\beta)}$ is the subalgebra of $\overline{aBr_m}(\beta)$ generated by all planar diagrams, i.e. diagrams that can be drawn without intersection of edges.
\item The \textit{affine Motzkin algebra} $\overline{aMo_m(\beta)}$ is the subalgebra of $\overline{aRoBr_m}(\beta)$ generated by all planar diagrams.
\item The \textit{affine planar rook algebra} $\overline{aPRo_m(\beta)}=kaPRo_m$ is the subalgebra of $\overline{aRo_m(\beta)}$ generated by all planar diagrams with no non-contractible loops.
\end{enumerate}
As before, {a component of a diagram means a component of the corresponding set partition}.
\end{Definition}

\begin{Remark}\label{planar no loop}
To avoid confusion, as above, planar in this paper means no intersecting edges ({where non-contractible loops count as edges}). We use the word planar to stay close to the traditional terminology in diagram algebras. {In particular, the planarity condition implies there can be no non-contractible loops for the planar affine diagram algebras unless there are no through lines, i.e. edges connecting a top vertex with a bottom one.}
\end{Remark}

The affine diagram algebras defined above should be seen as the `affinization' of the diagram algebras tabulated in \cite[Figure (1E.2)]{KhSiTu-monoidal-cryptography} and reproduced below in \autoref{fig: table of algebras}. Here affinization means `taking diagrams of a particular type and putting them on a cylinder or annulus': as above, the infinite horizontal strip in our definition is just a cylinder cut open; the edges that go across the boundaries between different translates of the fundamental rectangle are edges that wind around the cylinder, and the non-contractible loops are loops around the cylinder. We refer to the diagram algebras in the left column as the planar diagram algebras, and the ones in the right column as symmetric diagram algebras. We have omitted the planar partition algebra in our definition above because $pPa_{m}(\beta)$ is isomorphic to $TL_{2m}(\beta)$. Note that $pS_m$ is the trivial group, and its affinization is $\Z$ which can be seen by identifying
\begin{gather}\label{Z idenfification}
1\leftrightsquigarrow
\begin{tikzpicture}[anchorbase,scale=1]
\fill[pattern=north east lines,pattern color=black!30]
(-0.5,-0.05) rectangle (0.5,1.05);
\draw[usual] (0,0) to (-0.5,0.5);
\draw[usual] (0.5,0.5) to (0,1);
\draw[usual] (-0.5,0) to (0.5,0);
\draw[usual] (-0.5,1) to (0.5,1);
\draw[very thick,dashed,draw=gray] (-0.5,-0.05) to (-0.5,1.05);
\draw[very thick,dashed,draw=gray] (0.5,-0.05) to (0.5,1.05);
\end{tikzpicture},\quad
-1\leftrightsquigarrow
\begin{tikzpicture}[anchorbase,xscale=-1]
\fill[pattern=north east lines,pattern color=black!30]
(-0.5,-0.05) rectangle (0.5,1.05);
\draw[usual] (0,0) to (-0.5,0.5);
\draw[usual] (0.5,0.5) to (0,1);
\draw[usual] (-0.5,0) to (0.5,0);
\draw[usual] (-0.5,1) to (0.5,1);
\draw[very thick,dashed,draw=gray] (-0.5,-0.05) to (-0.5,1.05);
\draw[very thick,dashed,draw=gray] (0.5,-0.05) to (0.5,1.05);
\end{tikzpicture}
,
\end{gather}
and the affinization of $S_m$ is the wreath product $\Z \wr S_m$, with the same extra diagrammatic generators, which is sometimes called the extended affine symmetric group.

\begin{figure}[H]
\begin{tabular}{c|c|c||c|c|c}
\arrayrulecolor{tomato}
Symbol & Diagrams & Name: ${}_{-}$ algebra
& Symbol & Diagrams & Name: ${}_{-}$ algebra
\\
\hline
\hline
$pPa_m$ & \begin{tikzpicture}[anchorbase]
\draw[usual] (0.5,0) to[out=90,in=180] (1.25,0.45) to[out=0,in=90] (2,0);
\draw[usual] (0.5,0) to[out=90,in=180] (1,0.35) to[out=0,in=90] (1.5,0);
\draw[usual] (0.5,1) to[out=270,in=180] (1,0.55) to[out=0,in=270] (1.5,1);
\draw[usual] (1.5,1) to[out=270,in=180] (2,0.55) to[out=0,in=270] (2.5,1);
\draw[usual] (0,0) to (0,1);
\draw[usual] (2.5,0) to (2.5,1);
\draw[usual,dot] (1,0) to (1,0.2);
\draw[usual,dot] (1,1) to (1,0.8);
\draw[usual,dot] (2,1) to (2,0.8);
\end{tikzpicture} & Planar partition
& $Pa_m$ & \begin{tikzpicture}[anchorbase]
\draw[usual] (0.5,0) to[out=90,in=180] (1.25,0.45) to[out=0,in=90] (2,0);
\draw[usual] (0.5,0) to[out=90,in=180] (1,0.35) to[out=0,in=90] (1.5,0);
\draw[usual] (0,1) to[out=270,in=180] (0.75,0.55) to[out=0,in=270] (1.5,1);
\draw[usual] (1.5,1) to[out=270,in=180] (2,0.55) to[out=0,in=270] (2.5,1);
\draw[usual] (0,0) to (0.5,1);
\draw[usual] (1,0) to (1,1);
\draw[usual] (2.5,0) to (2.5,1);
\draw[usual,dot] (2,1) to (2,0.8);
\end{tikzpicture} & Partition
\\
\hline
$Mo_m$ & \begin{tikzpicture}[anchorbase]
\draw[usual] (0.5,0) to[out=90,in=180] (1.25,0.5) to[out=0,in=90] (2,0);
\draw[usual] (1,0) to[out=90,in=180] (1.25,0.25) to[out=0,in=90] (1.5,0);
\draw[usual] (2,1) to[out=270,in=180] (2.25,0.75) to[out=0,in=270] (2.5,1);
\draw[usual] (0,0) to (1,1);
\draw[usual,dot] (2.5,0) to (2.5,0.2);
\draw[usual,dot] (0,1) to (0,0.8);
\draw[usual,dot] (0.5,1) to (0.5,0.8);
\draw[usual,dot] (1.5,1) to (1.5,0.8);
\end{tikzpicture} & Motzkin
& $RoBr_m$ & \begin{tikzpicture}[anchorbase]
\draw[usual] (1,0) to[out=90,in=180] (1.25,0.25) to[out=0,in=90] (1.5,0);
\draw[usual] (1,1) to[out=270,in=180] (1.75,0.55) to[out=0,in=270] (2.5,1);
\draw[usual] (0,0) to (0.5,1);
\draw[usual] (2.5,0) to (2,1);
\draw[usual,dot] (0.5,0) to (0.5,0.2);
\draw[usual,dot] (2,0) to (2,0.2);
\draw[usual,dot] (0,1) to (0,0.8);
\draw[usual,dot] (1.5,1) to (1.5,0.8);
\end{tikzpicture} & Rook Brauer
\\
\hline
$TL_m$ & \begin{tikzpicture}[anchorbase]
\draw[usual] (0.5,0) to[out=90,in=180] (1.25,0.5) to[out=0,in=90] (2,0);
\draw[usual] (1,0) to[out=90,in=180] (1.25,0.25) to[out=0,in=90] (1.5,0);
\draw[usual] (0,1) to[out=270,in=180] (0.25,0.75) to[out=0,in=270] (0.5,1);
\draw[usual] (2,1) to[out=270,in=180] (2.25,0.75) to[out=0,in=270] (2.5,1);
\draw[usual] (0,0) to (1,1);
\draw[usual] (2.5,0) to (1.5,1);
\end{tikzpicture} & Temperley--Lieb
& $Br_m$ & \begin{tikzpicture}[anchorbase]
\draw[usual] (0.5,0) to[out=90,in=180] (1.25,0.45) to[out=0,in=90] (2,0);
\draw[usual] (1,0) to[out=90,in=180] (1.25,0.25) to[out=0,in=90] (1.5,0);
\draw[usual] (0,1) to[out=270,in=180] (0.75,0.55) to[out=0,in=270] (1.5,1);
\draw[usual] (1,1) to[out=270,in=180] (1.75,0.55) to[out=0,in=270] (2.5,1);
\draw[usual] (0,0) to (0.5,1);
\draw[usual] (2.5,0) to (2,1);
\end{tikzpicture} & Brauer
\\
\hline
$pRo_m$ & \begin{tikzpicture}[anchorbase]
\draw[usual] (0,0) to (0.5,1);
\draw[usual] (0.5,0) to (1,1);
\draw[usual] (2,0) to (1.5,1);
\draw[usual] (2.5,0) to (2.5,1);
\draw[usual,dot] (1,0) to (1,0.2);
\draw[usual,dot] (1.5,0) to (1.5,0.2);
\draw[usual,dot] (0,1) to (0,0.8);
\draw[usual,dot] (2,1) to (2,0.8);
\end{tikzpicture} & Planar rook
& $Ro_m$ & \begin{tikzpicture}[anchorbase]
\draw[usual] (0,0) to (1,1);
\draw[usual] (0.5,0) to (0,1);
\draw[usual] (2,0) to (2,1);
\draw[usual] (2.5,0) to (0.5,1);
\draw[usual,dot] (1,0) to (1,0.2);
\draw[usual,dot] (1.5,0) to (1.5,0.2);
\draw[usual,dot] (1.5,1) to (1.5,0.8);
\draw[usual,dot] (2.5,1) to (2.5,0.8);
\end{tikzpicture} & Rook
\\
\hline
$pS_m$ & \begin{tikzpicture}[anchorbase]
\draw[usual] (0,0) to (0,1);
\draw[usual] (0.5,0) to (0.5,1);
\draw[usual] (1,0) to (1,1);
\draw[usual] (1.5,0) to (1.5,1);
\draw[usual] (2,0) to (2,1);
\draw[usual] (2.5,0) to (2.5,1);
\end{tikzpicture} & Planar symmetric
& $S_m$ & \begin{tikzpicture}[anchorbase]
\draw[usual] (0,0) to (1,1);
\draw[usual] (0.5,0) to (0,1);
\draw[usual] (1,0) to (1.5,1);
\draw[usual] (1.5,0) to (2.5,1);
\draw[usual] (2,0) to (2,1);
\draw[usual] (2.5,0) to (0.5,1);
\end{tikzpicture} & Symmetric
\\
\end{tabular}
.
\caption{A table of diagram algebras.}
\label{fig: table of algebras}
\end{figure}

\begin{Remark}
The elements of the affine rook monoid can also be viewed as partial extended affine permutations, which have been studied in e.g.\cite{chmutov2018matrix}.    
\end{Remark}

\begin{Definition}
If $\overline{a\mathcal{X}_m(\beta)}$ is one of the affine diagram algebras defined above, we denote by $a\mathcal{X}_m(\beta,\alpha),\alpha\in k$ the corresponding \textbf{reduced affine} diagram algebra, which is generated by $m$-diagrams without non-contractible loops, and with multiplication defined so that if non-contractible loops appear they are removed with factor $\alpha$. 
\end{Definition}

Say an affine diagram is \textit{ordinary} if all vertices in the fundamental rectangle are connected only to vertices within it, and if it has no non-contractible loops. The ordinary diagrams form a subalgebra in $\overline{a\mathcal{X}_m}(\beta)$ and $a\mathcal{X}_m(\beta,\alpha)$ isomorphic to the non-affine diagram algebra $\mathcal{X}_m(\beta)$. 

\begin{Remark}
Rather than removing all internal components with the scalar $\beta$, it is possible to assign different parameters for the removal of different types of components, e.g. depending on the fundamental group of {the} component viewed as a graph, or depending on winding (cf. \autoref{oriented brauer section}).
\end{Remark}

\begin{Definition}
We also define the \textbf{periodic} diagram algebras  $\overline{p\mathcal{X}_m}(\beta)$ as follows: If $\mathcal{X}_m$ is one of the ordinary diagram algebras (e.g. Temperley--Lieb, Brauer, etc.), define $\overline{p\mathcal{X}_m}(\beta)$ by adding additional generators to the presentation of $\mathcal{X}$ and interpreting the old relations mod $m$. {The definition depends on a chosen presentation for each of the ordinary diagram algebras; the precise versions used here are defined below.}
\end{Definition}

For example, the Temperley--Lieb algebra $TL_m(\beta)$ has generators $e_1,\dots,e_{m-1}$ and relations $e_i^2=\beta e_i$,
$e_ie_{i\pm 1}e_i=e_i$,
$e_ie_j=e_je_i$, $\lvert i-j\rvert >1$. Then $\overline{pTL_m}(\beta)$ has a presentation obtained from the one above by adjoining a generator $e_0$ and interpreting all the relations mod $m$ (so that, for example, $e_0e_{m-1}e_0=e_{0}$). We will explain what presentations for each algebra we have in mind in \autoref{pres_alg} below. 

\begin{Definition} 
As before, we also define the reduced versions $p\mathcal{X}_m(\beta,\alpha)$, where we remove the non-contractible loops with the scalar factor $\alpha$.
\end{Definition}

The (reduced) periodic diagram algebras are subalgebras of the corresponding {(reduced) affine} diagram algebras. In the Temperley--Lieb case, the embedding maps the generator $e_i, 1\le i\le m-1$ to the corresponding ordinary diagram in $\overline{aTL_m}(\beta)$, and $e_0$ to the cup-cap generator that crosses over an edge of the fundamental domain, shown in \autoref{fig: e0 generator} below. The embedding for other algebras is similar.  

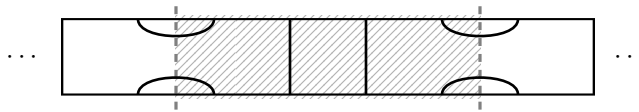
\begin{figure}[H]
\begin{center}
\begin{tikzpicture}[anchorbase]
\def\period{4} 
\begin{scope}

\clip (-2,-0.2) rectangle (5,1.2);

\fill[pattern=north east lines,pattern color=black!30]
(-0.5,-0.05) rectangle (3.5,1.05);
\draw[very thick,dashed,draw=gray] (-0.5,-0.2) -- (-0.5,1.2);
\draw[very thick,dashed,draw=gray] ( 3.5,-0.2) -- ( 3.5,1.2);
\draw[thick] (-2.5,0) -- (5.5,0);
\draw[thick] (-2.5,1) -- (5.5,1);

\foreach \k in {-1,0,1} {
\begin{scope}[xshift={\k*\period cm}]

\draw[usual] (1,0) -- (1,1);
\draw[usual] (2,0) -- (2,1);

\draw[usual]
(3,1) .. controls (3,0.7) and (4,0.7) .. (4,1);

\draw[usual]
(3,0) .. controls (3,0.3) and (4,0.3) .. (4,0);
\end{scope}
}

\end{scope}

\node at (-2.5,0.5) {$\cdots$};
\node at ( 5.5,0.5) {$\cdots$};
\end{tikzpicture}
\end{center}
\caption{The image of $e_0$ under the embedding $\overline{pTL_4}(\beta)\hookrightarrow \overline{aTL_4}(\beta)$.}
\label{fig: e0 generator}
\end{figure}

\begin{Notation}We write $\overline{a\mathcal{X}_m},a\mathcal{X}_m,\overline{p\mathcal{X}_m},$ and $p\mathcal{X}_m$ for the underlying monoids, so that when $\beta=1$ the algebras are the monoid algebras for these monoids.
\end{Notation}

\begin{Remark} 
For $\overline{a\mathcal{X}_m}$ and $a\mathcal{X}_m$ the group of units is either $\Z$ (planar case) or $\Z \wr S_m$ (symmetric case). For $\overline{p\mathcal{X}_m}$ and $p\mathcal{X}_m$ the group of units is either trivial (planar case) or the affine symmetric group $\widetilde{S}_m$ (symmetric case). 
\end{Remark}

The affine algebras are generated by the corresponding periodic algebras plus an additional invertible element $\tau$ of infinite order, representing a `global twist'. In the symmetric case this is equivalent to introducing $m$ winding elements $t_1,\dots, t_m$.

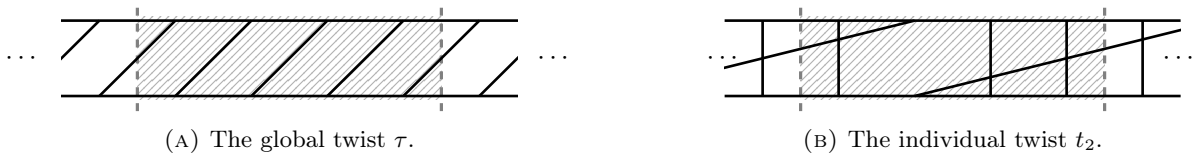
\begin{figure}[H]
\centering

\begin{subfigure}{0.48\textwidth}
\centering
\begin{tikzpicture}[anchorbase]
\def\period{4} 
\begin{scope}
\clip (-1.5,-0.2) rectangle (4.5,1.2);
\fill[pattern=north east lines,pattern color=black!30]
(-0.5,-0.05) rectangle (3.5,1.05);
\draw[very thick,dashed,draw=gray] (-0.5,-0.2) -- (-0.5,1.2);
\draw[very thick,dashed,draw=gray] ( 3.5,-0.2) -- ( 3.5,1.2);
\draw[usual] (-2.5,0) -- (5.5,0);
\draw[usual] (-2.5,1) -- (5.5,1);
\foreach \k in {-1,0,1} {
\begin{scope}[xshift={\k*\period cm}]
\foreach \j in {0,...,3} {
\draw[usual] (\j,0) -- (\j+1,1);
}
\end{scope}
}
\end{scope}
\node at (-2,0.5) {$\cdots$};
\node at ( 5,0.5) {$\cdots$};
\end{tikzpicture}
\caption{The global twist $\tau$.}
\label{fig:global-twist}
\end{subfigure}
\hfill
\begin{subfigure}{0.48\textwidth}
\centering
\begin{tikzpicture}[anchorbase]

\def\period{4}

\begin{scope}
\clip (-1.5,-0.2) rectangle (4.5,1.2);

\fill[pattern=north east lines,pattern color=black!30]
(-0.5,-0.05) rectangle (3.5,1.05);

\draw[very thick,dashed,draw=gray] (-0.5,-0.2) -- (-0.5,1.2);
\draw[very thick,dashed,draw=gray] ( 3.5,-0.2) -- ( 3.5,1.2);

\draw[usual] (-1.5,0) -- (4.5,0);
\draw[usual] (-1.5,1) -- (4.5,1);

\foreach \x in {-1,0,2,3,4} {
\draw[usual] (\x,0) -- (\x,1);
}

\draw[usual] (1,1) -- (-3,0);
\draw[usual] (5,1) -- (1,0);

\end{scope}

\node at (-1.5,0.5) {$\cdots$};
\node at ( 4.5,0.5) {$\cdots$};

\end{tikzpicture}
\caption{The individual twist $t_2$.}
\label{fig:individual-twist}
\end{subfigure}

\caption{The global twist $\tau$ and the individual twist $t_2$.}
\label{fig:twists}
\end{figure}

\subsubsection{Categories}\label{cat def section}
We can also define categories {corresponding} to these algebras. Our definition of an affine partition $m$-diagram can be generalized to define an affine partition $(s,t)$-diagram, where $s,t\in \Z_{\ge 0}$. Namely, we take the $\Z$-translation of a `fundamental rectangle' with $s$ vertices on the bottom edge and $t$ vertices on the top edge, and as before require that there be at most one edge between two vertices, and that the diagram be invariant under left or right translation.

\begin{Definition}
Let $k$ be a field and $\beta\in k$. Define the \textit{affine partition category} $\overline{{aPa}}^k(\beta)$ to be the $k$-linear category whose objects are the nonnegative integers, and define the morphism space from $s$ to $t$ to be the free $k$-vector space on the affine partition $(s,t)$-diagrams. Define the composition of morphisms by stacking and removing closed components with $\beta$, then the endomorphism spaces of $\overline{{aPa}}^k(\beta)$ are the algebras $\overline{aPa}_m(\beta)$. (Our convention is that $a: m\to n$ and $b:k\to m$ are composed to give $a\circ b: k\to n$.) 

Similarly, we define corresponding categories for any of the (reduced) affine or periodic algebras above. We also define the set-theoretic versions of the categories, where {the Hom sets are just the set of diagrams (rather than linear combinations thereof), and} endomorphism {sets} are monoids rather than $k$-algebras.  
\end{Definition}

\subsection{Relations to affine diagram algebras in the literature}

{The following is a far from exhaustive list} of affine diagram categories, algebras and monoids in the literature and how they relate to the above constructions.

\begin{Remark}
There are ways to represent the affine diagrams other than on an infinite horizontal strip or a cylinder. For example, they may be visualized using braids with a flagpole, see e.g. \cite{MR1604452,MR1824165,tubbenhauer2022handlebody}.
\end{Remark}

\begin{Remark}
As a word of warning, sometimes affine means attaching Jucys--Murphy-type elements or similar objects, often represented by dots, see e.g. \cite{MR1398116}. That is kind of the same, cf. \cite{tubbenhauer2022handlebody}, but most often not in a monoid-way as there are scalars and sums in relations.
We will not use any of these in this paper.
\end{Remark}

\subsubsection{Affine Temperley--Lieb}

Our $\overline{aTL_m}(\beta)$ is the same as the `extended' version of affine Temperley--Lieb algebra studied in e.g. \cite{graham1998representation} (denoted $\mathbf{T}^a(m)$) and \cite{green2023representations} (denoted $D_m$), while our $aTL_m(\beta)$ is known in the literature as the affine or periodic Temperley--Lieb algebra, and is a quotient of the affine Hecke algebra of type $\widetilde{A}_{m-1}$. We refer to \cite[Table 1]{langlois2023uncoiled} for a table comparing the terminology used in various papers. 

In \cite[\S 3]{langlois2023uncoiled} finite-dimensional quotients of $\overline{aTL}$ and $aTL$ are defined, called \textit{uncoiled} algebras; these impose finiteness conditions on winding in addition to getting rid of non-contractible loops. When $m$ is odd, the Jones annular algebra corresponds to the uncoiled $\overline{aTL}$ with unwinding parameter 1. (We note that \cite{langlois2023uncoiled} calls $\overline{aTL}$ affine and $aTL$ periodic.)  

The affinization of monoidal categories is studied in \cite{mousaaid2021affinization}. In particular, the affine Temperley--Lieb category of \cite{graham1998representation}, which corresponds to our $\overline{aTL}^k(\beta)$, can be viewed as the affinization of the ordinary Temperley--Lieb category.  

\subsubsection{Affine Brauer}\label{oriented brauer section}

The relationship between $aBr_m(\beta,\alpha)$ and the `colored Brauer algebra' $\widehat{D}_m(\delta_0,\delta_1,\dots)$ considered in \cite{goodman2004affine} and \cite{cui2014affine} is as follows.

\cite[\S 4.2]{goodman2004affine} defined a colored Brauer $m$-diagram to be a Brauer diagram (cf. \autoref{fig: table of algebras}) where each edge is labeled with an integer. Let the $2m$ vertices of the Brauer diagram be ordered by $1<\dots < m < \overline{m} <\dots <\overline{1}$, where $1,\dots,m$ are the top vertices and $\overline{1},\dots,\overline{m}$ are the bottom ones. Two colored Brauer $m$-diagrams can be stacked to produce a new diagram, where the label of each new edge is obtained by traversing the composite strand and summing the labels of the \textit{oriented strands} encountered (that is, if an edge of label $r$ goes from a lower index to a higher one then it contributes $r$, and $-r$ if the other way around). The algebra $\widehat{D}_m(\delta_0,\delta_1,\dots)$ is then the $k$-algebra with a basis given by colored Brauer $m$-diagrams, and multiplication given by stacking and then removing closed components with the scalars $\delta_i\in k, i\in \Z_{\ge 0}$, where $i$ is the absolute value of the sum of labels as one traverses a given closed component.

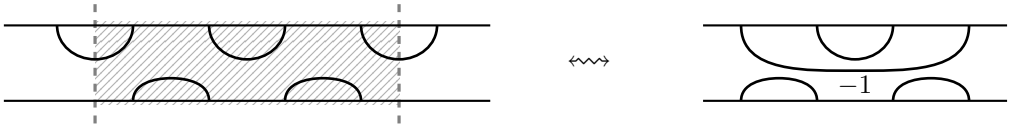
\begin{figure}[H]
\begin{center}
\begin{tikzpicture}[anchorbase]
\begin{scope}[xshift=-6cm]
\begin{scope}
\clip (-0.7,-0.3) rectangle (5.7,1.3);
\fill[pattern=north east lines,pattern color=black!30]
(0.5,-0.05) rectangle (4.5,1.05);
\draw[very thick,dashed,gray] (0.5,-0.3) -- (0.5,1.3);
\draw[very thick,dashed,gray] (4.5,-0.3) -- (4.5,1.3);
\draw[thick] (-0.7,0) -- (5.7,0);
\draw[thick] (-0.7,1) -- (5.7,1);
\draw[usual] (0,1)
to[out=270,in=180] (0.5,0.55)
to[out=0,in=270] (1,1);
\draw[usual] (2,1)
to[out=270,in=180] (2.5,0.55)
to[out=0,in=270] (3,1);
\draw[usual] (4,1)
to[out=270,in=180] (4.5,0.55)
to[out=0,in=270] (5,1);
\draw[usual] (1,0)
to[out=90,in=180] (1.5,0.30)
to[out=0,in=90] (2,0);
\draw[usual] (3,0)
to[out=90,in=180] (3.5,0.30)
to[out=0,in=90] (4,0);
\end{scope}
\end{scope}
\node at (1,0.5) {\Large$\leftrightsquigarrow$};
\begin{scope}[xshift=3cm]
\draw[thick] (-0.5,0) -- (3.5,0);
\draw[thick] (-0.5,1) -- (3.5,1);
\draw[usual] (0,1)
to[out=270,in=180] (1.5,0.40)
to[out=0,in=270]  (3,1);
\draw[usual] (1,1)
to[out=270,in=180] (1.5,0.55)
to[out=0,in=270]  (2,1);
\node[below] at (1.5,0.45) {$-1$};
\draw[usual] (0,0)
to[out=90,in=180] (0.5,0.30)
to[out=0,in=90]  (1,0);
\draw[usual] (2,0)
to[out=90,in=180] (2.5,0.30)
to[out=0,in=90]  (3,0);
\end{scope}
\end{tikzpicture}
\end{center}
\caption{The affine diagram on the left corresponds to the colored diagram on the right. The labels which are $0$ are not shown.}
\label{fig: colored iso}
\end{figure}

Analogously, we may define such algebras (call them \textit{colored diagram algebras}) associated with the rook and rook-Brauer algebras, where we assign an integer label to every edge; if we set $\delta_i=\beta$ for all $i$, then this algebra is isomorphic to the reduced affine algebra $a\mathcal{X}_m(\beta, \alpha=\beta)$, where $\mathcal{X}_m$ is one of $Ro_m, RoBr_m,$ or $Br_m$. To see this, note that the data of a reduced affine diagram consist exactly of an ordinary diagram together with an integer for each edge, specifying how many boundaries between translates of the fundamental rectangle has been crossed (to the left or right). See \autoref{fig: colored iso} for an example illustrating this isomorphism: note that the label of $-1$ corresponds to a translation of one rectangle to the left; this is negative because in the corresponding diagram on the left the edge is negatively oriented. Under this isomorphism, the `individual twist' $t_i$ is just the identity diagram with the $i$-th strand given label $1$. We also note that planarity is not preserved by this isomorphism, so the correct version of planar affine diagram algebras is the one in \autoref{affine diagram def}.

{It is more complicated to define a colored partition algebra which is isomorphic to $aPa_m(\beta)$. In this case, there is no canonical set of edges, and we should label each pair of vertices $(u,v)$ in a component (possibly $u=v$) of a partition diagram by a nonempty set of integers $D(u,v)$, thought of as recording their relative positions; we need this to be a set because in an affine partition diagram a vertex may be connected to more than one translate of a vertex in the fundamental rectangle. Again, each label in the set is to be understood in an oriented sense, $D(u,v)=-D(v,u):=\{-x\mid x\in D(v,u)\}$. Moreover, if $u,v,w$ are vertices in a component, we require $D(u,v)+D(v,w):=\{x+y\mid x\in D(u,v), y\in D(v,w) \}=D(u,w)$. Note that the construction forces $D(u,u)$ to be a subgroup of $\Z$; moreover, this subgroup depends only on the component $u$ is in. We now define composition similarly to before: first stack the underlying diagrams $P, Q$ and remove each closed component (with a factor of $\beta$). Suppose the vertices $x$ in $P$ and $y$ in $Q$ are in the same component of the resulting diagram; for a \textit{chain from $x$ to $y$}, i.e. a sequence $\gamma= (x=x_0, B_1, x_1, B_2,\dots, B_r, x_r=y)$ where $B_i$ is a block of one of the two original diagrams and $x_{i-1},x_i\in B_i$, we obtain a well-defined set $D(\gamma)=D(x_0,x_1)+\dots+D(x_{r-1},x_r)$. We then define $D(x,y)$ in the resulting diagram to be the union over all such chains. It is not hard to check that this is well-defined and gives an algebra isomorphic to $aPa_m(\beta)$.}

{
We note that our construction here is similar to the notion of \textit{group partition categories} in \cite{likeng2021group}, where each vertex is labelled by a group element modulo uniform translation of labels in each component; however, the composition defined in that paper uses a different convention which requires the compatibility of labels.}

\begin{Notation}
We write $a\mathcal{X}_m(\beta)$ for $a\mathcal{X}_m(\beta, \alpha=\beta)$.
\end{Notation}

\subsubsection{Others}

As far as we can tell, for all the other affine diagram monoids there is no reference in the same sense as above.

\subsection{Presentations for algebras}\label{pres_alg}

In this section we give the presentations for all the periodic monoids, as well as for $\overline{aTL_m}(\beta)$, $aBr_m(\beta),aRo_m(\beta)$, $aRoBr_m(\beta)$ and ${aPa_m}(\beta)$. For $\overline{aMo_m}(\beta)$ and ${aPRo_m}(\beta)$ we give a generating set. These results will be used in \autoref{fin gen nonunit}. {Throughout we assume $m\ge 3$. When we say a relation with the condition $|i-j|>1$ is \textit{interpreted mod $m$}, this means the condition becomes $j\not\equiv i,i\pm 1 \pmod m$.}
\subsubsection{Temperley--Lieb}
The periodic Temperley--Lieb algebra $\overline{pTL_m}(\beta)$, see \cite[Definition 2.2.1]{green2023representations}, is the algebra with generators $e_0,\dots, e_{m-1}$ and the following relations (where $0\le i\le m-1$ and indices are interpreted mod $m$): 
\begin{enumerate}
\item $e_i^2=\beta e_i$,
\item $e_ie_{i\pm 1}e_i=e_i$,
\item $e_ie_j=e_je_i$, $\lvert i-j\rvert >1$. 
\end{enumerate}
The affine Temperley--Lieb algebra $\overline{aTL}_m(\beta)$ (assume $m\ge 3$) has an additional invertible generator $\tau$ and additional relations:
\begin{enumerate}[resume]
\item $\tau e_i\tau^{-1}=e_{i+1}$ (indices taken mod $m$).
\item $(\tau e_1)^{(m-1)}=\tau^m (\tau e_1)$. This is equivalent to the relation $e_{i+2}e_{i+3}...e_{i}=\tau^2 e_i$.
\end{enumerate}
This presentation may be found in \cite[Proposition 2.3.7]{green2023representations}. 

Say an affine Temperley--Lieb diagram is \textit{even} if the number of times it intersects with the line $i+1/2$ for any $i$ is an even number (the `number of intersections' is defined to be the minimal possible after perturbing the edges). The algebra $\overline{pTL_m}(\beta)$ is known to have the following diagrammatic interpretation:

\begin{Proposition}
$\overline{pTL_m}(\beta)$ is the subalgebra of $\overline{aTL_m}(\beta)$ generated by all even diagrams, excluding {nonzero} even powers of $\tau$.    
\end{Proposition}
\begin{proof}
See \cite[Proposition 2.2.3]{green2023representations} or \cite[Proposition 2.2]{langlois2023uncoiled}. 
\end{proof}

\subsubsection{Planar rook}\label{planar rook pres sec}
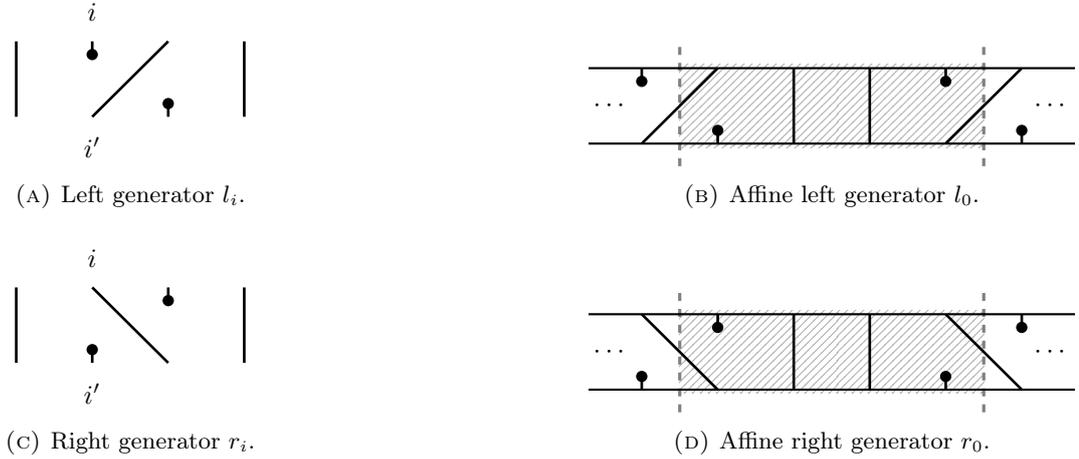
\begin{figure}[H]
\centering


\begin{subfigure}{0.45\textwidth}
\centering
\begin{tikzpicture}[anchorbase]

\draw[usual] (0,0) -- (0,1);
\draw[usual] (3,0) -- (3,1);

\draw[usual] (2,1) -- (1,0);

\draw[usual,dot] (1,1) to (1,0.8);  
\draw[usual,dot] (2,0) to (2,0.2);  

\node[above] at (1,1.15) {$i$};
\node[below] at (1,-0.15) {$i'$};

\end{tikzpicture}
\caption{Left generator $l_i$.}
\end{subfigure}
\hfill
\begin{subfigure}{0.45\textwidth}
\centering
\begin{tikzpicture}[anchorbase]

\def\period{4}

\begin{scope}
\clip (-0.7,-0.3) rectangle (5.7,1.3);

\fill[pattern=north east lines,pattern color=black!30]
(0.5,-0.05) rectangle (4.5,1.05);

\draw[very thick,dashed,gray] (0.5,-0.3) -- (0.5,1.3);
\draw[very thick,dashed,gray] (4.5,-0.3) -- (4.5,1.3);

\draw[thick] (-0.7,0) -- (5.7,0);
\draw[thick] (-0.7,1) -- (5.7,1);

\foreach \k in {-1,0,1} {
\pgfmathsetmacro{\S}{3 + 4*\k}

\draw[usual] (\S+0,0) -- (\S+0,1);
\draw[usual] (\S+3,0) -- (\S+3,1);

\draw[usual] (\S+2,1) -- (\S+1,0);

\draw[usual,dot] (\S+1,1) to (\S+1,0.8); 
\draw[usual,dot] (\S+2,0) to (\S+2,0.2); 
}

\end{scope}

\node at (-0.4,0.5) {$\cdots$};
\node at (5.4,0.5) {$\cdots$};

\end{tikzpicture}
\caption{Affine left generator $l_0$.}
\end{subfigure}

\vspace{1em}


\begin{subfigure}{0.45\textwidth}
\centering
\begin{tikzpicture}[anchorbase]

\draw[usual] (0,0) -- (0,1);
\draw[usual] (3,0) -- (3,1);

\draw[usual] (1,1) -- (2,0);

\node[above] at (1,1.15) {$i$};
\node[below] at (1,-0.15) {$i'$};

\draw[usual,dot] (2,1) to (2,0.8); 
\draw[usual,dot] (1,0) to (1,0.2); 

\end{tikzpicture}
\caption{Right generator $r_i$.}
\end{subfigure}
\hfill
\begin{subfigure}{0.45\textwidth}
\centering
\begin{tikzpicture}[anchorbase]

\def\period{4}

\begin{scope}
\clip (-0.7,-0.3) rectangle (5.7,1.3);

\fill[pattern=north east lines,pattern color=black!30]
(0.5,-0.05) rectangle (4.5,1.05);

\draw[very thick,dashed,gray] (0.5,-0.3) -- (0.5,1.3);
\draw[very thick,dashed,gray] (4.5,-0.3) -- (4.5,1.3);

\draw[thick] (-0.7,0) -- (5.7,0);
\draw[thick] (-0.7,1) -- (5.7,1);

\foreach \k in {-1,0,1} {
\pgfmathsetmacro{\S}{3 + 4*\k}

\draw[usual] (\S+0,0) -- (\S+0,1);
\draw[usual] (\S+3,0) -- (\S+3,1);

\draw[usual] (\S+1,1) -- (\S+2,0);

\draw[usual,dot] (\S+2,1) to (\S+2,0.8); 
\draw[usual,dot] (\S+1,0) to (\S+1,0.2); 
}

\end{scope}

\node at (-0.4,0.5) {$\cdots$};
\node at (5.4,0.5) {$\cdots$};

\end{tikzpicture}
\caption{Affine right generator $r_0$.}
\end{subfigure}

\caption{Finite and affine versions of the left and right generators $l_i$ and $r_i$.}
\label{fig: left right gen}
\end{figure}

The presentation of the planar rook monoid $PRo_m$ can be found in \cite[\S 3.3]{herbig2006planar}; the presentation for the algebra $PRo_m(\beta)$ is similar. The generators are $l_i, r_i$, $1\le i \le m-1$, with relations:
\begin{enumerate}
\item $l_i^3 = \beta l_i^2=\beta r_i^2=r_i^3$, 
\item $r_i r_{i+1} r_i = r_i r_{i+1} = r_{i+1} r_i r_{i+1}$, $l_i l_{i+1} l_i = l_{i+1} l_i = l_{i+1} l_i l_{i+1}$, 
\item $r_i l_i r_i = r_i$, $l_i r_i l_i = l_i$, 
\item $r_{i+1} l_i r_i = r_{i+1} l_i$, $l_{i-1} r_i l_i = l_{i-1} r_i$,  
\item $l_i r_i l_{i+1} = r_i l_{i+1}$, $r_i l_i r_{i-1} = l_i r_{i-1}$, 
\item $r_i l_i = l_{i+1} r_{i+1}$, 
\item if $|i-j|\geq 2$, then $r_i l_j = l_j r_i$, $r_i r_j = r_j r_i$, $l_i l_j = l_j l_i$. 
\end{enumerate}
We thus define $pPRo_m(\beta)$ to be generated by $l_i,r_i$, $0\le i\le m-1$, with relations from above but interpreted mod $m$. The generators $l_i, r_i$ are illustrated in \autoref{fig: left right gen}.

\begin{Proposition}\label{planar rook gen}
The affine planar rook algebra $aPRo_m(\beta)$ is generated by $l_1, r_1, \tau^{\pm 1}$, where $\tau$ is the diagram illustrated in \autoref{fig:twists}. Moreover, the following relations (with indices interpreted mod $m$) hold in addition to the $pPRo_m$ relations above: \begin{enumerate}[resume]
\item $\tau l_i \tau^{-1}=l_{i+1}, (\tau l_1)^m=\tau^m(\tau l_1)$,
\item $\tau r_i\tau^{-1}=r_{i+1}, (\tau r_1)^m=\tau^m(\tau r_1)$.
\end{enumerate}
{Note that we only claim that $\{l_1, r_1, \tau^{\pm 1}\}$ is a generating set, and that the above relations hold; we do not claim that the above gives a presentation.}
\end{Proposition}
\begin{proof}
The relations are easily verified by drawing diagrams. In particular, they imply that $l_1,r_1, \tau^{\pm 1}$ generate all the $l_i$, $r_i$. We now show that any diagram $D$ in $aPRo_m$ can be produced from these. If $D$ has no isolated components in the fundamental rectangle, then planarity implies $D$ is a power of $\tau$. Otherwise, if there are no edges emanating from vertices in the rectangle which connect two vertices of {at least $m$ units apart}, then the diagram is just the translate of an ordinary planar rook diagram and can be generated from the $l_i, r_i$ as in \cite[\S 3.2, Lemma 1]{herbig2006planar}. However, if there are longer edges, planarity means that by multiplying by a suitable power of $\tau$ we can always reduce to the previous case, so we are done. 
\end{proof}

\subsubsection{Motzkin}
The presentation of the Motzkin algebra $Mo_m(\beta)$ can be found in \cite[Theorem 4.1]{posner2013presentation}. It has generators $e_i,l_i,r_i, 1\le i \le m-1$ and the following relations, in addition to the planar rook relations above in \autoref{planar rook pres sec}: 
\begin{enumerate}
\item If $|i-j|\geq 2$, then $e_i r_j = r_j e_i$, $e_i l_j = l_j e_i$, $e_i e_j = e_j e_i$, 
\item $e_i^2 = \beta e_i$, $e_i e_{i\pm 1} e_i = e_i$, 
\item $e_i r_i = e_i l_i$, $r_i e_i = l_ie_i $, 
\item $e_i l_{i+1} = e_{i} e_{i+1} r_i$, $r_{i+1} e_i = l_i e_{i+1} e_i$, 
\item $l_i l_{i+1} e_i = e_{i+1} l_i l_{i+1}$, $e_ir_{i+1}r_i = r_{i+1} r_i e_{i+1}$, 
\item $e_i r_i e_i = \beta e_i$. 
\end{enumerate}
We define $\overline{pMo_m}(\beta)$ to be the algebra with generators $e_i,l_i,r_i, 0\le i\le m-1$ and all relations of the Motzkin algebra interpreted mod $m$.

\begin{Proposition}\label{motzkin gen}
The algebra $\overline{aMo_m} (\beta)$ is generated by $e_1, l_1, r_1, \tau^{\pm 1}$. Moreover, the following relations (with indices interpreted mod $m$) hold in addition to the $\overline{pMo_m}(\beta)$ relations: \begin{enumerate}[resume]
\item $\tau l_i \tau^{-1}=l_{i+1}$, $(\tau l_1)^m=\tau^m(\tau l_1)$,
\item $\tau r_i\tau^{-1}=r_{i+1}$, $(\tau r_1)^m=\tau^m(\tau r_1)$.
\item $\tau e_i\tau^{-1}=e_{i+1},(\tau e_1)^{(m-1)}=\tau^m (\tau e_1)$.
\end{enumerate}   
{Note that we only claim that $\{e_1, l_1, r_1, \tau^{\pm 1}\}$ is a generating set, and that the above relations hold; we do not claim that the above gives a presentation.}
\end{Proposition}
\begin{proof}
The proof is similar to that of \autoref{planar rook gen}.  
\end{proof}

\subsubsection{Brauer}
The Brauer algebra $Br_m(\beta)$ is generated by $s_i, e_i, 1\le i \le m-1$, with the following relations (see e.g. \cite[Lemma 2.10]{lehrer2012brauercategoryinvarianttheory}):
\begin{enumerate}
\item $s_i^2=1$, $s_is_j=s_js_i$, $\lvert i-j\rvert >1$, $s_is_{i+1}s_i=s_{i+1}s_is_{i+1}$,
\item $e_i^2 = \beta e_i$, $e_i e_{i\pm 1} e_i = e_i$,
\item If $\lvert i-j\rvert >1$, $e_i e_j = e_j e_i$, $s_i e_j = e_j s_i$,
\item $s_i s_{i\pm 1} e_i = e_{i\pm 1} e_i$, $e_i s_{i\pm 1} s_i = e_i e_{i\pm 1}$, $s_i e_i = e_i s_i = e_i$.
\end{enumerate}
We thus define $\overline{pBr_m}(\beta)$ to be the algebra generated by $s_i,e_i,0\le i\le m-1$ with the above relations interpreted mod $m$. 

\begin{Proposition}\label{Brauer pres}
The algebra ${aBr_m}(\beta)$ is generated by $e_i,s_i, 1\le i \le m-1$, and $t_i^{\pm 1}, 1\le i \le m$, subject to the relations of ${Br_m}(\beta)$ together with the following, giving a list of defining relations:\begin{enumerate}[resume]
  \item {$t_it_i^{-1}=t_i^{-1}t_i=1$,}
\item $s_it_j=t_js_i$, {$j\notin \{i,i+1\}$},
\item $e_it_j=t_je_i$, {$j\notin \{i,i+1\}$},
\item $t_it_j=t_jt_i$,
\item $s_it_i=t_{i+1}s_i$,
\item $e_it_it_{i+1}=e_i=t_it_{i+1}e_i$,
\item $e_it_i^ae_i=\beta e_i$ for $a\in \Z$.
\end{enumerate}    
\end{Proposition}
The element $t_i$ should be viewed as increasing the winding number of the $i$-th strand by one.
\begin{proof}
  This presentation first appeared in \cite[\S 3.2]{cui2014affine} with hints toward a proof. To prove it, it is helpful to recall the interpretation that a reduced affine Brauer diagram is a Brauer diagram with edges labelled {by integers}, cf. \autoref{oriented brauer section}. {Any word in the generators given above} can be put in the `normal form' ${\beta^d}(t_1^{i_1}\dots t_{m}^{i_{m}})D(t_1^{j_1}\dots t_{m}^{j_{m}})$, where $D$ is an ordinary Brauer diagram and $i_k, j_k \in \Z$, {and $d$ is the number of closed components removed in the composition (of the underlying Brauer diagrams). That the relation $e_it_i^ae_i=\beta e_i$ is enough to remove all loops can be proven by imitating the proof that $e_i^2=\beta e_i$ is enough in the ordinary Brauer case, see \cite[Lemma A.7]{lehrer2012brauercategoryinvarianttheory}.}  Since the presentation of Brauer algebras is included in the proposed presentation, it suffices to show that if a diagram has two such normal forms with the same $D$, then they can be transformed into one another via the additional relations. Now, the relations $s_it_j=t_js_i$ if $j\notin\{i,i+1\}$ and $s_it_i=t_{i+1}s_i$ imply that commuting the $t_i$'s through the $s_j$'s always transports the winding to the correct strand, and the relations $e_it_j=t_je_i$, $j\notin\{i,i+1\}$ make sure that the $t_i$'s behave correctly with the cups/caps, namely: they commute with cups that are far away, if the two vertices of a cup/cap have the same amount of winding then that cancels out (because the winding of a component is \textit{relative} to its two vertices). Using these relations, one can move all the windings attached to through lines to one side, and the windings on the cups/caps must also agree, so the two normal forms are equivalent. {This completes the proof.}   
\end{proof}

\subsubsection{Rook}
By \cite[Theorem 44.3]{lipscomb1996symmetric}, the rook algebra $Ro_m(\beta)$ is generated by $e, s_i, 1\le i\le m-1$, and relations
\begin{enumerate}
\item $s_i^2=1$, $s_is_j=s_js_i$, $\lvert i-j\rvert >1$, $s_is_{i+1}s_i=s_{i+1}s_is_{i+1}$,
\item $e^2=\beta e,(es_1)^3=\beta (es_1)^2$,
\item $es_i=s_ie$, $i>1$,
\end{enumerate}
where $e$ represents the identity diagram with the first through line removed, and $es_1$ corresponds to $l_1$ from before.
We define $pRo_m(\beta)$ to be the algebra with generators $e, s_i, 0\le i\le m-1$, and the relations above (interpreted mod $m$, {where $es_i=s_ie$ now applies to $i\neq 0, 1$}) {together with the additional relation $(es_0)^3=\beta(es_0)^2$}.

\begin{Proposition}The algebra $aRo_m(\beta)$ has a presentation with generators $e, s_i, 1\le i\le m-1$ and $t_i^{\pm 1}, 1\le i\le m$, and the relations for $Ro_m$ together with
\begin{enumerate}[resume]
\item {$t_it_i^{-1}=t_i^{-1}t_i=1$,}
\item $s_it_j=t_js_i$, $j\notin \{i,i+1\}$,
\item $t_it_j=t_jt_i$,
\item $s_it_i=t_{i+1}s_i$,
\item $t_1e=et_1=e$, and $t_ie=et_i$, $i>1$.
\end{enumerate}    
\end{Proposition}
\begin{proof}
The proof is similar to that of \autoref{Brauer pres}.     
\end{proof}

\subsubsection{Rook-Brauer}\label{rook brauer section}
By \cite[Theorem 5.1]{mazorchuk2005presentations}, the rook-Brauer algebra $RoBr_m(\beta)$ is generated by $e_i, s_i, 1\le i \le m-1$ and Brauer relations together with generators $p_i, 1\le i \le m$ (isolated components in the $i$-th position and vertical lines elsewhere), with additional relations:
\begin{enumerate}
\item {$p_i^2=\beta p_i$, $p_ip_j=p_jp_i, i\neq j$,}
\item {$s_ip_i=p_{i+1}s_i$, $p_is_ip_i=p_ip_{i+1}$, $s_ip_j=p_js_i, j\notin \{i,i+1\}$,}  
\item $e_ip_j=p_je_i$, $j\notin \{i,i+1\}$,
\item $e_ip_i=e_ip_{i+1}, {e_ip_ip_{i+1}=\beta e_ip_i}, p_ie_i=p_{i+1}e_i, {p_{i}p_{i+1}e_i=\beta p_ie_i}$,
\item $e_ip_ie_i=\beta e_i$, $p_ie_ip_i=p_ip_{i+1}.$
\end{enumerate}
We thus define $\overline{pRoBr_m}(\beta)$ to be generated by $e_i, s_i, p_i, 0\le i \le m-1$, subject to the relations of $RoBr_m(\beta)$ interpreted mod $m$.

\begin{Proposition}\label{rook brauer pres}
The algebra $aRoBr_m(\beta)$ has a presentation with generators $e_i, s_i, 1\le i \le m-1$ and $p_i, t_i^{\pm 1}, 1\le i\le m$, with the relations of $aBr_m(\beta)$ in \autoref{Brauer pres}, the relations of $RoBr_m(\beta)$, and the additional relations $t_ip_i=p_it_i=p_i$ {and $t_ip_j=p_jt_i, i\neq j$}.   
\end{Proposition}
\begin{proof}
The proof is similar to that of \autoref{Brauer pres}. The additional relations make sure that isolated components do not have windings.    
\end{proof}

\subsubsection{Partition}
By \cite[Theorem 1.11]{halverson2005partition} ({see also \cite[\S 6.3]{east2011generators}}), the partition algebra $Pa_m(\beta)$ is generated by $s_i, 1\le i\le m-1$ and $p_j$, {$j\in\{1,3/2,\dots, m\}$}, with the following relations:
\begin{enumerate}
\item $p_j^2 = \beta p_j, {p_{j+1/2}^2=p_{j+1/2},}$ {for integral $j$},
\item $p_jp_{j \pm 1/2}p_j = p_j,$ $p_ip_j = p_jp_i, |i-j| > 1/2,$
\item $s_i^2 = 1,$ $s_is_{i+1}s_i = s_{i+1}s_is_{i+1},$ $s_is_j = s_js_i, \quad \text{for } |i-j| > 1,$
\item $s_ip_ip_{i+1} = p_ip_{i+1}s_i = p_ip_{i+1},$ $s_ip_is_i = p_{i+1},$ 
\item $s_ip_{i+1/2} = p_{i+1/2}s_i = p_{i+1/2},$ $s_is_{i+1}p_{i+1/2}s_{i+1}s_i = p_{i+3/2},$
\item $s_ip_j = p_js_i, \quad \text{for } j \neq i-1/2, \, i, \, i+1/2, \, i+1, \, i+3/2.$
\end{enumerate}
(The $s_i$'s are the crossings, the $p_j$'s have isolated components in the $j$-th position, and $p_{i+1/2}$ has a square connecting the $i$-th and $(i+1)$-st vertices.)

We define $\overline{pPa_m}(\beta)$ to be the algebra with generators $s_i, p_i, p_{i+1/2}, 0\le i \le m-1$, subject to the relations above interpreted mod $m$. 

\begin{Proposition}

The algebra $aPa_m(\beta)$ has a presentation with generators
$s_i$, $p_r$ and $t_j^{\pm1}$, where
$
1\le i\le m-1,
r\in\{1/2,1,3/2,\dots,m\},
1\le j\le m,$
with the following relations in addition to those of $Pa_m(\beta)$:

\begin{enumerate}[resume]

\item

{$t_jt_j^{-1}=t_j^{-1}t_j=1$,
$t_it_j=t_jt_i$,}

\item

$s_it_i=t_{i+1}s_i$,
$s_it_j=t_js_i$ for $1\le i\le m-1$, $1\le j\le m$ and $j\notin \{i,i+1\}$,

\item

$p_rt_r=t_rp_r=p_r$ for $1\le r\le m$,
and
$p_rt_j=t_jp_r$ for $1\le r,j\le m$ and $j\neq r$,

\item

$p_{i+1/2}t_j=t_jp_{i+1/2}$ for $1\le i\le m-1$, $1\le j\le m$ and $j\notin \{i,i+1\}$,

\item

$p_{i+1/2}t_it_{i+1}
=t_it_{i+1}p_{i+1/2}$ for $1\le i\le m-1$.
\item {$p_ip_{i+1}p_{i+1/2}t_i^ap_{i+1/2}p_ip_{i+1}=\beta p_ip_{i+1}$ for $1\le i\le m-1$ and $a\in \Z$}.

\end{enumerate}

\end{Proposition}

\begin{proof}

The proof is similar to that of \autoref{Brauer pres}. The relations $p_{i+1/2}t_it_{i+1}=t_it_{i+1}p_{i+1/2}$ and $p_{i+1/2}t_j=t_jp_{i+1/2}$ for $j\neq i,i+1$ make sure that windings commute correctly with squares, and that again `only relative positions matter'.  Any connected component can be created with squares and point removals, and the relations keep track of how many points there are in each component: for example,  $p_{i+1/2}p_ip_{i+1}$ is a cup, and we have $t_it_{i+1}p_{i+1/2}p_ip_{i+1}=p_{i+1/2}t_it_{i+1}p_ip_{i+1}=p_{i+1/2}p_ip_{i+1}$, and so we recover the cup relations.     

\end{proof}

\begin{figure}[H]
\centering


\begin{subfigure}{0.45\textwidth}
\centering
\begin{tikzpicture}[anchorbase]

\draw[usual] (0,0) -- (0,1);
\draw[usual] (3,0) -- (3,1);

\draw[usual] (1,0) -- (2,1);
\draw[usual] (2,0) -- (1,1);

\node[above] at (1,1.15) {$i$};
\node[above] at (2,1.15) {$i+1$};
\node[below] at (1,-0.15) {$i'$};
\node[below] at (2,-0.15) {$(i+1)'$};

\end{tikzpicture}
\caption{Finite crossing $s_i$.}
\end{subfigure}
\hfill
\begin{subfigure}{0.45\textwidth}
\centering
\begin{tikzpicture}[anchorbase]

\begin{scope}
\clip (-0.7,-0.3) rectangle (5.7,1.3);

\fill[pattern=north east lines,pattern color=black!30]
(0.5,-0.05) rectangle (4.5,1.05);

\draw[very thick,dashed,gray] (0.5,-0.3) -- (0.5,1.3);
\draw[very thick,dashed,gray] (4.5,-0.3) -- (4.5,1.3);

\draw[thick] (-0.7,0) -- (5.7,0);
\draw[thick] (-0.7,1) -- (5.7,1);

\foreach \shift in {-1,3} {

\draw[usual] (\shift+0,0) -- (\shift+0,1);
\draw[usual] (\shift+3,0) -- (\shift+3,1);

\draw[usual] (\shift+1,0) -- (\shift+2,1);
\draw[usual] (\shift+2,0) -- (\shift+1,1);
}

\end{scope}

\node at (-0.4,0.5) {$\cdots$};
\node at ( 5.4,0.5) {$\cdots$};

\end{tikzpicture}
\caption{Affine crossing $s_0$.}
\end{subfigure}

\vspace{1em}


\begin{subfigure}{0.45\textwidth}
\centering
\begin{tikzpicture}[anchorbase]

\draw[usual] (0,0) -- (0,1);
\draw[usual] (3,0) -- (3,1);


\draw[usual] (1,0)
to[out=90,in=90] (2,0);

\draw[usual] (1,1)
to[out=270,in=270] (2,1);

\draw[usual] (1,0)
to[out=70,in=290] (1,1);

\draw[usual] (2,0)
to[out=110,in=250] (2,1);

\node[above] at (1,1.15) {$i$};
\node[above] at (2,1.15) {$i+1$};
\node[below] at (1,-0.15) {$i'$};
\node[below] at (2,-0.15) {$(i+1)'$};

\end{tikzpicture}
\caption{Finite rectangle $p_{i+1/2}$.}
\end{subfigure}
\hfill
\begin{subfigure}{0.45\textwidth}
\centering
\begin{tikzpicture}[anchorbase]

\begin{scope}
\clip (-0.7,-0.3) rectangle (5.7,1.3);

\fill[pattern=north east lines,pattern color=black!30]
(0.5,-0.05) rectangle (4.5,1.05);

\draw[very thick,dashed,gray] (0.5,-0.3) -- (0.5,1.3);
\draw[very thick,dashed,gray] (4.5,-0.3) -- (4.5,1.3);

\draw[thick] (-0.7,0) -- (5.7,0);
\draw[thick] (-0.7,1) -- (5.7,1);

\foreach \shift in {-1,3} {

\draw[usual] (\shift+0,0) -- (\shift+0,1);
\draw[usual] (\shift+3,0) -- (\shift+3,1);

\draw[usual] (\shift+1,0)
to[out=90,in=90] (\shift+2,0);

\draw[usual] (\shift+1,1)
to[out=270,in=270] (\shift+2,1);

\draw[usual] (\shift+1,0)
to[out=70,in=290] (\shift+1,1);

\draw[usual] (\shift+2,0)
to[out=110,in=250] (\shift+2,1);
}

\end{scope}

\node at (-0.4,0.5) {$\cdots$};
\node at ( 5.4,0.5) {$\cdots$};

\end{tikzpicture}
\caption{Affine rectangle $p_{i+1/2}$.}
\end{subfigure}

\caption{Finite and affine versions of the crossing $s_i$ and rectangle $p_{i+1/2}$.}
\end{figure}
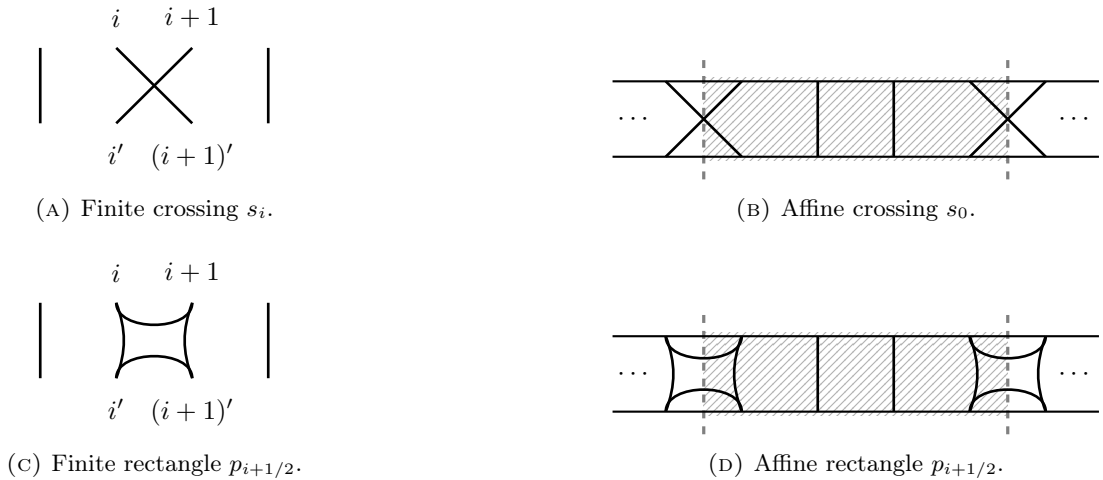

\subsection{Presentations for monoidal categories}

Geometrically, when the underlying non-affine categories are braided, then
\begin{gather*}
\begin{tikzpicture}[anchorbase,scale=1,yscale=1]
\fill[pattern=north east lines,pattern color=black!30]
(-0.5,-1.05) rectangle (1.5,1.05);
\draw[spinach!35,fill=spinach!35] (-0.5,-0.25) to (-0.5,0.25) to (1.5,0.25) to (1.5,-0.25) to (-0.5,-0.25);
\draw[usual] (0,-1) to (0,-0.25);
\draw[usual] (0,0.25) to (0,1);
\draw[usual] (1,-1) to (1,-0.25);
\draw[usual] (1,0.25) to (1,1);
\node at (0.5,-0.625) {$\dots$};
\node at (0.5,0.625) {$\dots$};
\draw[very thick,dashed,gray] (-0.5,-1.1) to (-0.5,1.1);
\draw[very thick,dashed,gray] (1.5,-1.1) to (1.5,1.1);
\node at (0.5,0) {$D$};
\end{tikzpicture}	
\otimes
\begin{tikzpicture}[anchorbase,scale=1,yscale=1]
\fill[pattern=north east lines,pattern color=black!30]
(-0.5,-1.05) rectangle (1.5,1.05);
\draw[spinach!35,fill=spinach!35] (-0.5,-0.25) to (-0.5,0.25) to (1.5,0.25) to (1.5,-0.25) to (-0.5,-0.25);
\draw[usual] (0,-1) to (0,-0.25);
\draw[usual] (0,0.25) to (0,1);
\draw[usual] (1,-1) to (1,-0.25);
\draw[usual] (1,0.25) to (1,1);
\node at (0.5,-0.625) {$\dots$};
\node at (0.5,0.625) {$\dots$};
\draw[very thick,dashed,gray] (-0.5,-1.1) to (-0.5,1.1);
\draw[very thick,dashed,gray] (1.5,-1.1) to (1.5,1.1);
\node at (0.5,0) {$E$};
\end{tikzpicture}
=
\begin{tikzpicture}[anchorbase,scale=1,yscale=1]
\fill[pattern=north east lines,pattern color=black!30]
(-0.5,-1.05) rectangle (3.5,1.05);
\draw[spinach!35,fill=spinach!35] (-0.5,0.25) to (3.5,0.25) to (3.5,0.75) to (-0.5,0.75) to (-0.5,0.25);
\draw[spinach!35,fill=spinach!35] (-0.5,-0.25) to (3.5,-0.25) to (3.5,-0.75) to (-0.5,-0.75) to (-0.5,-0.25);
\draw[usual] (0,-1) to (0,-0.75);
\draw[usual,crossline] (0,-0.25) to (0,1);
\draw[usual] (1,-1) to (1,-0.75);
\draw[usual,crossline] (1,-0.25) to (1,1);
\draw[usual] (2,-1) to (2,-0.85);
\draw[usual] (2,-0.15) to (2,0.25);
\draw[usual] (2,0.75) to (2,1);
\draw[usual] (3,-1) to (3,-0.85);
\draw[usual] (3,-0.15) to (3,0.25);
\draw[usual] (3,0.75) to (3,1);
\node at (0.5,-1) {$\dots$};
\node at (0.5,1) {$\dots$};
\node at (2.5,-1) {$\dots$};
\node at (2.5,1) {$\dots$};
\draw[very thick,dashed,gray] (-0.5,-1.1) to (-0.5,1.1);
\draw[very thick,dashed,gray] (3.5,-1.1) to (3.5,1.1);
\node at (0.5,-0.5) {$D$};
\node at (2.5,0.5) {$E$};
\end{tikzpicture}
\end{gather*}
defines a monoidal structure on affine versions.
This is easy when the braiding is a symmetry: The categories associated to the symmetric reduced affine diagram algebras have a (strict) monoidal structure: given two reduced affine symmetric diagrams, we think of them as ordinary symmetric diagrams with integer labels, and simply juxtapose them.
That the categories associated to affine Temperley--Lieb and affine Motzkin algebras are monoidal is more subtle and will be deduced from the above formula (via the helpful \cite[Theorem 2.6]{mousaaid2021affinization}), using the braiding on Temperley--Lieb and Motzkin categories. 

We now give monoidal presentations for the $k$-linear versions of these categories; the presentations for the set-theoretic analogs are similar.

\subsubsection{Partition}

\begin{Proposition}
As a strict $k$-linear monoidal category, $aPa^k(\beta)$  is generated by a single generating object 1, generating morphisms $\Xop: 2\to 2$ (crossing), $\Dop: 2\to 2$ (square), $\etop:1\to 0$ (top dot), $\ebot: 0\to 1$ (bottom dot), $\phi^+: 1\to 1$ (positive winding), $\phi^-: 1\to 1$ (negative winding), and the following defining relations (where we write $\iota_n$ for the identities and also $I$ for $\iota_1$.): 
\begin{enumerate}

\item $\Xop\circ \Xop= \iota_2, \etop \circ \ebot= \beta \iota_0$,

\item $\Dop\circ \Dop=\Dop, \Dop=\Dop\circ \Xop = \Xop\circ \Dop$, $\Dop\otimes I\circ I\otimes \Dop=I\otimes \Dop\circ \Dop\otimes I$

\item $(\Xop\otimes I)\circ (I\otimes \Xop)\circ (\Xop\otimes I)=(I\otimes \Xop)\circ(\Xop\otimes I)\circ (I\otimes \Xop)$,

\item $(\Xop\otimes I)\circ (I\otimes \Dop)\circ (\Xop\otimes I)=(I\otimes \Xop)\circ(\Dop\otimes I)\circ (I\otimes \Xop)$,

\item $\Xop\circ (I\otimes \ebot)=\ebot\otimes I, (I\otimes \etop)\circ \Xop=\etop \otimes I$,

\item $(I\otimes \etop)\circ \Dop\circ (I\otimes \ebot)=I$, $\Dop\circ (I\otimes \etop \otimes \ebot)\circ \Dop=\Dop$,

\end{enumerate}

\begin{enumerate}[resume]

\item $\phi^+ \circ \phi^- = I={\phi^-\circ \phi^+}$,

\item $(\phi^\pm \otimes I)\circ \Xop = \Xop\circ (I\otimes \phi^\pm)$, $(I\otimes \phi^{\pm})\circ \Xop= \Xop\circ (\phi^{\pm}\otimes I)$,

\item $\etop \circ \phi^{\pm}=\etop$, $\phi^{\pm}\circ \ebot = \ebot$,

\item $(\phi^{\pm}\otimes \phi^{\pm})\circ \Dop=\Dop\circ (\phi^{\pm}\otimes \phi^{\pm})$
\item
{$(\etop\otimes \etop)\circ \Dop\circ (\phi^a\otimes I)\circ \Dop\circ(\ebot\otimes \ebot)=\beta\iota_0$, $a\in \Z$.}
\end{enumerate}

\end{Proposition}
\begin{proof}
The generators $\Xop,\Dop,\ebot,\etop$ are known to generate the partition category with relations (a)-(f),  see \cite[Theorem 3.7]{east2023presentations}. That these together with the  additional generators $\phi^\pm$ and relations (g)-(k) generate $aPa^k(\beta)$ is proven just as for the algebra. 
\end{proof}

\subsubsection{Brauer}

\begin{Proposition}
As a strict $k$-linear monoidal category, $aBr^k(\beta)$ is generated by a single generating object $1$, generating morphisms 
$\Xop:2\to 2$, $\Cupop: 0\to 2$, $\Capop: 2\to 0$, $\phi^{\pm}:1\to 1$, and the following defining relations: 
\begin{enumerate}
\item $\Xop\circ \Xop= \iota_2$, $\Capop \circ \Cupop = \beta\iota_0$, $\Xop\circ \Cupop = \Cupop, \Capop \circ \Xop = \Capop$, 
\item $(\Xop\otimes I)\circ (I\otimes \Xop)\circ (\Xop\otimes I)=(I\otimes \Xop)\circ(\Xop\otimes I)\circ (I\otimes \Xop) $,
\item $(I\otimes \Capop) \circ (\Cupop \otimes I)=I=(\Capop \otimes I)\circ (I\otimes \Cupop)$,
\item $(\Xop\otimes I)\circ (I\otimes \Cupop) =  (I\otimes \Xop)\circ (\Cupop \otimes I)$, $(\Capop \otimes I)\circ (I\otimes \Xop)= (I\otimes \Capop)\circ (\Xop\otimes I)$.
\end{enumerate}
\begin{enumerate}[resume]
\item $\phi^+ \circ \phi^- = I={\phi^-\circ \phi^+}$,
\item $(\phi^\pm \otimes I)\circ \Xop = \Xop\circ (I\otimes \phi^\pm)$, $(I\otimes \phi^{\pm})\circ \Xop= \Xop\circ (\phi^{\pm}\otimes I)$,
\item $\Capop\circ (\phi^+\otimes \phi^+) = \Capop$, $ (\phi^+ \otimes \phi^+)\circ \Cupop=\Cupop$,
{\item $\Capop\circ(\phi^a\otimes I)\circ\Cupop
=\beta\iota_0$ for $a\in\Z$.}
\end{enumerate}
\end{Proposition}
\begin{proof}
Again, $\Xop,\Cupop,\Capop$ and the relations (a)-(d) generate the Brauer category, which is well-known, see e.g. \cite[Theorem 2.6]{lehrer2012brauercategoryinvarianttheory}, and also \cite[Theorem 3.2]{east2023presentations}. Now the result follows by a similar proof to that of \autoref{Brauer pres}.  
\end{proof}

\subsubsection{Rook}

\begin{Proposition}
As a strict $k$-linear monoidal category, $aRo^k(\beta)$ has a presentation as follows. It is generated by a single generating object 1, generating morphisms $\Xop:2\to 2,\etop: 1\to 0, \ebot: 0\to 1, \phi^\pm:1\to 1$, and all relations involving {only} these generators from the presentations of $aPa^k(\beta)$ and $aBr^k(\beta)$ above.
\end{Proposition}
\begin{proof}
Proven as for the algebra, see \cite[\S 5.2]{hu2019presentationsdiagramcategories} for presentation of the rook category. 
\end{proof}

\subsubsection{Rook-Brauer}
\begin{Proposition}
As a strict $k$-linear monoidal category, $aRoBr^k(\beta)$ has a presentation as follows. It is generated by a single generating object 1, generating morphisms $\Xop: 2\to 2$, $\Cupop:0\to 2$, $\Capop: 2\to 0$, $\etop: 1\to 0$, $\ebot: 0\to 1$, $\phi^{\pm}:1\to 1$, subject to the relations of $RoBr^k(\beta)$ together with all the additional relations involving $\phi^\pm$ from the presentations of $aRo^k(\beta)$ and $aBr^k(\beta)$ above.
\end{Proposition}
\begin{proof}
Proven as for the algebra, see \cite[\S 6.2]{hu2019presentationsdiagramcategories} for presentation of the rook-Brauer category.  
\end{proof}

\subsubsection{Temperley--Lieb and Motzkin}
We refer to \cite[\S 7.1, \S 7.2]{hu2019presentationsdiagramcategories} for the definitions and presentations of the Temperley--Lieb category $TL(\beta)$ and the Motzkin category $Mo(\beta)$. These categories have a braided structure, which can be used to turn their affinizations into strict monoidal categories.

\begin{Lemma}\label{braiding lemma} \
\begin{enumerate}
\item If $\beta$ is of the form $-q-q^{-1}$ where $q\in k^\times$ has a square root, then the Temperley--Lieb category $TL(\beta)$ admits a braiding $\eta$ determined by
\[\eta_{1,1}=q^{\frac12}\; \tikzvcenter{%
\draw[usual] (0,0) -- (0,0.8);
\draw[usual] (0.5,0) -- (0.5,0.8);
}+q^{-\frac12}\tikzvcenter{%
\draw[usual] (0,0.9)
to[out=270,in=180] (0.25,0.6)
to[out=0,in=270]  (0.5,0.9);
\draw[usual] (0,0.1)
to[out=90,in=180] (0.25,0.4)
to[out=0,in=90]   (0.5,0.1);
}\;,\]
with the inverse obtained by exchanging $q^{1/2}\leftrightarrow q^{-1/2}$ in $\eta_{1,1}$.

\item If $\beta$ is of the form $1-q-q^{-1}$, where $q\in k^\times$ has a square root, then the Motzkin category $Mo(\beta)$ admits a braiding $\sigma$ determined by
\[
\sigma_{1,1} =
q^{\frac12}\!\left(
\tikzvcenter{%
\draw[usual] (0,0) -- (0,0.8);
\draw[usual] (0.5,0) -- (0.5,0.8);
}
-
\tikzvcenter{%
\draw[usual] (0,0) -- (0,0.8);
\draw[usual,dot] (0.5,0) -- (0.5,0.2);
\draw[usual,dot] (0.5,0.8) -- (0.5,0.6);
}
-
\tikzvcenter{%
\draw[usual] (0.5,0) -- (0.5,0.8);
\draw[usual,dot] (0,0) -- (0,0.2);
\draw[usual,dot] (0,0.8) -- (0,0.6);
}
\right)
\;+\;
q^{-\frac12}\!\left(
\tikzvcenter{%
\draw[usual] (0,0.9)
to[out=270,in=180] (0.25,0.6)
to[out=0,in=270]  (0.5,0.9);
\draw[usual] (0,0.1)
to[out=90,in=180] (0.25,0.4)
to[out=0,in=90]   (0.5,0.1);
}
-
\tikzvcenter{%
\draw[usual] (0,0.9)
to[out=270,in=180] (0.25,0.6)
to[out=0,in=270]  (0.5,0.9);
\draw[usual,dot] (0,0.1)   -- (0,0.3);
\draw[usual,dot] (0.5,0.1) -- (0.5,0.3);
}
-
\tikzvcenter{%
\draw[usual] (0,0.1)
to[out=90,in=180] (0.25,0.4)
to[out=0,in=90]   (0.5,0.1);
\draw[usual,dot] (0,0.9)   -- (0,0.7);
\draw[usual,dot] (0.5,0.9) -- (0.5,0.7);
}
\right)
\;+\;
\tikzvcenter{%
\draw[usual] (0.5,0.8) -- (0,0);
\draw[usual] (0,0.8)   -- (0.5,0);
}
\;+\;
(q^{-\tfrac12}+q^{\tfrac12}-1)\,
\tikzvcenter{%
\draw[usual,dot] (0,0.1)   -- (0,0.3);
\draw[usual,dot] (0.5,0.1) -- (0.5,0.3);
\draw[usual,dot] (0,0.9)   -- (0,0.7);
\draw[usual,dot] (0.5,0.9) -- (0.5,0.7);
}\;,
\]
with the inverse obtained by exchanging $q^{1/2}\leftrightarrow q^{-1/2}$ in $\sigma_{1,1}$. 
\end{enumerate}
\end{Lemma}
\begin{proof}
The braiding on the Temperley--Lieb category is well-known, see \cite[\S 8.4]{mousaaid2021affinization} and \cite[Proposition 2.6]{belletete2018fusion}. The braiding on the \textit{dilute Temperley--Lieb category} $dTL(\beta)$ defined in \cite{belletete2018fusion} can be found in \cite[Proposition 5.1]{belletete2018fusion}. We have that $dTL(\beta-1)\cong Mo(\beta)$, and the formula for $\sigma_{1,1}$ above is obtained by translating dilute Temperley--Lieb diagrams to Motzkin diagrams via the algebra isomorphism \cite[\S 9.2]{nurcombe2024algebraicaspectslatticemodels} $\phi: dTL_m(\beta-1)\xrightarrow{\sim} Mo_m(\beta)$ that takes a diagram $D$ to \[\sum_{s=0}^{s(D)} (-1)^s \sum_{D' \in D(s)} D',\]
where $s(D)$ is the number of strings in $D$, and $D(s)$ is the set of diagrams obtained from $D$ by removing $s$ strings.
\end{proof}

\begin{Proposition}\label{motz cat pres}
{Suppose $\beta$ satisfies the hypotheses of \autoref{braiding lemma}.} The affine Temperley--Lieb category $\overline{aTL}^k(\beta)$ and the affine Motzkin category $\overline{aMo}^k(\beta)$ are strict monoidal categories, which are obtained from $TL(\beta)$ and $Mo(\beta)$ by adjoining invertible morphisms $\zeta_i: i\to i$ for each object $i\in \Z_{\ge 0}$, subject to the relations that \begin{enumerate}
\item $\alpha_{i,j}(\zeta_i\otimes \iota_j)=(\iota_j\otimes \zeta_i) \alpha_{j,i}^{-1}$,
\item $\zeta_{i\otimes j}=\zeta_i \otimes \zeta_j$,
\item $\zeta_j f=f\zeta_i$ 
\end{enumerate}  
for all objects $i,j$ and morphisms $f:i\to j$, where as before we write $\iota$ for the identity morphism, and $\alpha$ is either $\eta$ (for Temperley--Lieb) or $\sigma$ (for Motzkin) from \autoref{braiding lemma}. These give presentations.
\end{Proposition}
\begin{proof}
This follows from \cite[Theorem 2.6]{mousaaid2021affinization}. 
\end{proof}
{In the presentation above, the generators $\zeta_i$ correspond to winding $i$ strands around the cylinder, see the diagrams in \cite[Definition 2.1]{mousaaid2021affinization}.}

The presentation for $\overline{aTL}^k(\beta)$ may be simplified, see \cite[\S 8.4]{mousaaid2021affinization} for details.

\section{Representation theory}\label{rep section}

In this section we classify the simple modules of the planar affine and periodic diagram algebras, and those of the affine and periodic rook algebra. We also study the representation theory of the other symmetric affine or periodic diagram algebras by considering their finite-dimensional quotients, whose simple modules are described. Our main tool is sandwich cellularity.

All modules are left modules and are assumed to be finite-dimensional. {We adopt the convention that the symmetric group on 0 letters, $S_0$, is the trivial group. Similarly, the (extended) affine symmetric group on $0$ letters is understood to be the trivial group.}

\subsection{The (extended) affine symmetric groups}\label{semidrect backgroun} We begin by briefly discussing the simple representations over $\C$ of the `groups of units' of the symmetric affine diagram algebras, which will be used to construct simple modules for the affine and periodic rook algebras over $\C$. 
 Recall that the affine symmetric group $\widetilde{S}_m$ is the group with generators $s_0,\dots, s_{m-1}$ and relations
\begin{enumerate}
\item $s_i^2 = 1$ for all $0 \leq i \le m-1$, 
\item $s_i s_j = s_j s_i$ whenever $|i-j| \geq 2$, 
\item $s_i s_{i+1} s_i = s_{i+1} s_i s_{i+1}$ for all $0 \leq i \le m-1$, 

\end{enumerate}
where indices are taken modulo $m$. The picture to keep in mind is (here $m=2$):
\begin{gather*}
s_0\leftrightsquigarrow
\begin{tikzpicture}[anchorbase,scale=1]
\fill[pattern=north east lines,pattern color=black!30]
(-0.25,-0.05) rectangle (0.75,0.55);
\draw[usual] (0,0) to (-0.5,0.5);
\draw[usual] (-0.5,0) to (0,0.5);
\draw[usual] (0.5,0) to (1,0.5);
\draw[usual] (1,0) to (0.5,0.5);
\draw[very thick,dashed,draw=gray] (-0.25,-0.05) to (-0.25,0.55);
\draw[very thick,dashed,draw=gray] (0.75,-0.05) to (0.75,0.55);
\end{tikzpicture}
,\quad
s_1\leftrightsquigarrow
\begin{tikzpicture}[anchorbase,scale=1]
\fill[pattern=north east lines,pattern color=black!30]
(-0.25,-0.05) rectangle (0.75,0.55);
\draw[usual] (0,0) to (0.5,0.5);
\draw[usual] (0.5,0) to (0,0.5);
\draw[very thick,dashed,draw=gray] (-0.25,-0.05) to (-0.25,0.55);
\draw[very thick,dashed,draw=gray] (0.75,-0.05) to (0.75,0.55);
\end{tikzpicture}.
\end{gather*}
Algebraically, it can be viewed as the semidirect product $\Lambda \rtimes S_m$, where $\Lambda$ is the type A root lattice $\Lambda = \{(a_1,\dots, a_m)\in \Z^m\mid \sum a_i=0\}\cong \Z^{m-1}$ viewed as a subgroup of $\Z^m$, and $S_m$ acts on it by permutation of coordinates. If we adjoin to $\widetilde{S}_m$ a `global twist' element $\tau$ depicted as before as
\begin{gather*}
\tau\leftrightsquigarrow
\begin{tikzpicture}[anchorbase]
\def\period{4} 
\begin{scope}
\clip (-1.5,-0.2) rectangle (4.5,1.2);
\fill[pattern=north east lines,pattern color=black!30]
(-0.5,-0.05) rectangle (3.5,1.05);
\draw[very thick,dashed,draw=gray] (-0.5,-0.2) -- (-0.5,1.2);
\draw[very thick,dashed,draw=gray] ( 3.5,-0.2) -- ( 3.5,1.2);
\draw[usual] (-2.5,0) -- (5.5,0);
\draw[usual] (-2.5,1) -- (5.5,1);
\foreach \k in {-1,0,1} {
\begin{scope}[xshift={\k*\period cm}]
\foreach \j in {0,...,3} {
\draw[usual] (\j,0) -- (\j+1,1);
}
\end{scope}
}
\end{scope}
\node at (-2,0.5) {$\cdots$};
\node at ( 5,0.5) {$\cdots$};
\end{tikzpicture}
,
\end{gather*}
then we obtain the extended affine symmetric group $\Z\wr S_m$.

Recall that in characteristic zero, the simple representations of semidirect products by a finite index abelian subgroup are easily described. Suppose $A$ is a finite index abelian subgroup of $G, G=A\rtimes H$, then $G$ acts on the character group $ X:=\Hom(A,\C^\times)$ of $A$ by \[s\chi(a)=\chi(s^{-1}as), s\in G, \chi \in X, a\in A.\]
Let $(\chi_i)_{i \in X/H}$ be representatives for the orbits of $H$ in $X$. For each $i$ let $H_i$ be the stabilizer subgroup of $\chi_i$ in $H$, and let $G_i:=A\cdot H_i$ be the corresponding \textit{inertial subgroup} in $G$. We can extend $\chi_i$ to a character $\tilde{\chi}_i$ of $G_i$ by setting \[\tilde{\chi}_i(ah)=\chi_i(a), a\in A,h \in H_i.\]
 Now if $\rho$ is a simple representation of $H_i$, by composing with the projection $G_i \to H_i$ we get a simple representation $\tilde{\rho}$ of $G_i$. Denote by $\theta_{i,\rho}$ the induction of the tensor product $\tilde{\chi}_i\otimes \tilde{\rho}$ to $G$.

\begin{Proposition}\label{semidirect thm serre}
Let $G=A\rtimes H$ where $A$ is abelian and $H$ is finite. All simple $\C G$-modules are of the form $\theta_{i,\rho}$ and these are pairwise nonisomorphic.
\end{Proposition}

\begin{proof}
Well-known, e.g. the proof given in \cite[\S 8.2]{serre1977linear} for the case where $G$ is finite also works in the general case where $G$ is potentially infinite and $A$ is of finite index. (The proof only uses Mackey's criterion and the fact that the restriction of a simple module to a normal subgroup is semisimple: both facts remain true when $G$ is infinite and the subgroups are of finite index, the latter by Clifford's theorem.)
\end{proof}

\begin{Lemma}\label{rep of wreath}
Let $G=\Z \wr S_m$. Let $Par$ denote the set of all integer partitions $\mu$. The simple $\C G$-modules are all constructed as above and are parameterized by \[\Gamma_m:=\{f: \C^\times \to Par : \sum_{\alpha \in \C^\times}\lvert f(\alpha)\rvert = m\}.\]

Moreover, if for $f\in \Gamma_m$ we have $f(\alpha)=\mu_{\alpha}, n_\alpha:=\lvert\mu_{\alpha}\rvert$, and $L_f$ is the simple module indexed by $f$, then \[\dim L_f= \frac{m!}{\prod_\alpha n_\alpha!}\prod_\alpha \dim S^{\mu_\alpha},\] where $S^{\mu_\alpha}$ is the Specht module corresponding to the partition $\mu_\alpha$. 
\end{Lemma}

\begin{proof}
The simple $\C G$-modules are parameterized by pairs $(\chi,\rho)$ where $\chi$ {represents a $S_m$-orbit of} characters of $\Z^m$ and $\rho$ is a simple module of its inertial subgroup. The inertial subgroup of $\chi$ is a young subgroup {$\prod_\alpha S_{n_\alpha}, \sum_{n_\alpha}=m$ where $n_\alpha$ counts the number of basis vectors of $\Z^m$ which are sent to $\alpha$ under $\chi$.} A simple module of this group {is determined by a partition $\mu_\alpha\vdash n_\alpha$} for each $\alpha$ occurring in $\chi$. It follows that $\Gamma_m$ parameterizes the simple $\C G$-modules, and the formula for the dimension follows from the construction of simple modules. 
\end{proof}

\begin{Remark}
For the affine symmetric group $\widetilde{S}_m \cong \Z^{m-1}\rtimes S_m$ the parameterization is somewhat more complicated. The characters of $\Z^{m-1}$ can be viewed as the set of $m$-tuples with value in {$\C ^\times$}, quotiented by the relation that identifies scalar multiples. For generic characters $\chi$ the inertial subgroup is still a young subgroup, but the inertial subgroup of say $\chi = (1,\omega,\omega^2, 1,\omega,\omega^2)$, where $\omega^3=1$, is $C_2\wr C_3$. Nevertheless, complex simple representations of $\widetilde{S}_m$ still all have the form described in \autoref{semidirect thm serre}{, as $\widetilde{S}_m \cong \Z^{m-1}\rtimes S_m$}.  
\end{Remark}

\subsection{Sandwich cellular algebras}

Sandwich cellular algebras (see \cite[Definition 2.2]{tubbenhauer2022handlebody}, \cite[Definition 2A.3]{Tubbenhauer_2024}) are a generalization of the cellular algebras of \cite{lehrer1996cellular} and the affine cellular algebras of \cite{koenig2012affine}. We will only need the special case where the sandwich cellular algebra is also \textit{involutive}.

\begin{Definition} \label{sandwich def}
An \textbf{involutive sandwich cellular algebra} is an associative $k$-algebra $A$ with a tuple $(\Lambda, \mathcal{T},(H_\lambda,B_\lambda, {(\rule{0.5em}{0.4pt})^{\star_\lambda}}),C, (\rule{0.5em}{0.4pt}
)^\star)$, called the involutive sandwich cell datum, where
\begin{itemize}
\item $(\Lambda,\le)$ is a poset,
\item $\mathcal{T}=\bigcup_{\lambda\in \Lambda} \mathcal{T}(\lambda)$ is a collection of finite sets indexed by $\Lambda$, called top sets,
\item For $\lambda \in \Lambda$, $H_\lambda$ is a $k$-algebra with basis $B_\lambda$, and 
\item $C:\coprod_{\lambda \in \Lambda} \mathcal{T}(\lambda)\times B_\lambda\times \mathcal{T}(\lambda)\to A; (T,m,B)\mapsto c_{T,m,B}^\lambda$ is an injective map, 
\item $(\rule{0.5em}{0.4pt}
)^\star$ is an anti-automorphism of $A${, and for $\lambda\in \Lambda, (\rule{0.5em}{0.4pt})^{\star_\lambda}:B_\lambda\to B_\lambda$ is an order 2 bijection.} 
\end{itemize}
such that the following properties are satisfied:
\begin{enumerate}
\item[(AC1)] The set $B_A=\{c^\lambda_{T,m,B}: \lambda \in \Lambda, m\in B_\lambda, T,B \in \mathcal{T}(\lambda)\}$ is a basis of $A$.
\item[(AC2)] For $x\in A$, there exist scalars $r_{T,U,n}^x\in k$ independent of $B$ or $m$, such that \begin{equation}
xc_{T,m,B}^\lambda \equiv \sum_{U\in \mathcal{T}(\lambda),n\in B_{\lambda}} r_{T,U,n}^xc_{U,n,B}^\lambda \quad (\bmod A^{>\lambda}),
\end{equation}
where $A^{>\lambda}$ is the $k$-submodule of $A$ spanned by $\{c_{U,n,B}^\mu : \mu \in \Lambda, \mu > \lambda, U, B\in \mathcal{T}(\mu)\}$. A similar {relation} is required to hold for right multiplication. 
\item[(AC3)] There is an $A$-$H_\lambda$-bimodule $\Delta(\lambda)$ {which is free as a right $H_\lambda$-module}, an $H_{\lambda}$-$A$-bimodule $\nabla(\lambda)$ {which is free as a left $H_{\lambda}$-module}, and an $A$-bimodule isomorphism \begin{equation} \label{Alambda}
A_{\lambda}:=A^{\ge \lambda}/A^{>\lambda}\cong \Delta(\lambda)\otimes_{H_\lambda}\nabla(\lambda),
\end{equation}
where $A^{\ge \lambda} $ is the $k$-submodule of $A$ defined similarly to $A^{>\lambda}$. We call $\Delta(\lambda)$ and $\nabla(\lambda)$ the left and right cell modules, and $A_{\lambda}$ the cell algebra.
\item[(AC4)] We have $(c_{T,m,B}^\lambda)^\star\equiv c_{B,m^{\star_\lambda},T}^\lambda \quad (\bmod A^{>\lambda}).$  
\end{enumerate}
We also say an algebra $A$ is \textbf{involutive $\infty$ sandwich cellular} if it satisfies all of the above but the top and bottom sets are not necessarily finite.
\end{Definition}

Note that an algebra may be involutive ($\infty$) sandwich cellular in various ways. When we say below that an algebra is (not) involutive sandwich cellular we mean with respect to the cell data specified. 

The (reduced) affine and (reduced) periodic diagram algebras are all involutive $\infty$ sandwich cellular, and are all involutive sandwich cellular except those associated to Brauer, Rook-Brauer, and partition algebras. We now explain the ($\infty$) sandwich cellular structures. The poset $\Lambda$ for these algebras {is} either $\Lambda_1:=\{\lambda \in \Z_{\ge 0}, \lambda \le m\}$ or $\Lambda_2:=\{\lambda \in \Z_{\ge 0}, \lambda \le m, \lambda \equiv m \bmod 2\}$, with the \textit{reverse} of the usual linear order. The ($\infty$) sandwich cell data for the algebras are summarized in \autoref{tab:sandwich}. 

{In all cases, let $D(\lambda,m)$ denote the set of $(\lambda,m)$-diagrams (as defined in \autoref{cat def section}) with $\lambda$ through lines and no non-contractible loops. Let $G_\lambda$ denote the group of invertible $(\lambda,\lambda)$-diagrams. There is a natural free right action of $G_\lambda$ on $D(\lambda,m)$ by precomposition, and we take the top sets $\mathcal{T}(\lambda)$ to be a chosen set of representatives for the orbits $D(\lambda,m)/G_\lambda$. For example, in the planar case we take $\mathcal{T}(\lambda)$ to be the set of $(\lambda,m)$-diagrams with $\lambda$ through lines which have no `winding' in the through lines. The involution $(\rule{0.5em}{0.4pt}
)^\star$ is defined by reflecting the diagram across a horizontal line through the middle of the diagram.}

As is usual for diagram algebras, an affine (say partition) diagram $D$ with $t$ non-contractible loops has {(once the representatives $\mathcal{T}(\lambda)$ are chosen)} a unique decomposition into top, middle and bottom components, $D=D_u\circ D_\lambda \circ D_v$, where $D_u\in \mathcal{T}(\lambda)$ is a $(\lambda,m)$-diagram, $D_\lambda$ is a $(\lambda,\lambda)$-diagram with $t$ non-contractible loops and $\lambda$ through lines, $D_v$ is a $(m,\lambda)$-diagram with $D_v^\star\in \mathcal{T}(\lambda)$, and $\lambda$ is the minimal nonnegative integer for which such a decomposition exists. Extend the involution $(\rule{0.5em}{0.4pt}
)^\star$ to the category $\overline{aPa}^k(\beta)$. Let $H_\lambda$ be the subalgebra with a basis given by $B_\lambda$, the set of affine partition $(\lambda, \lambda)$-diagrams with $\lambda$ through lines {(possibly with non-contractible loops)}. {The order 2 bijection $(\rule{0.5em}{0.4pt}
)^{\star_\lambda}$ is again given by reflection.} We obtain a map $C:\coprod_{\lambda \in \Lambda_1}\mathcal{T}(\lambda) \times B_\lambda \times \mathcal{T}(\lambda) \to \overline{aPa_m}(\beta)$ which sends $(D_u,D_{\lambda}, D_{b})$ to the diagram $D=D_u \circ D_\lambda \circ (D_b)^\star$, whose image is a basis of the algebra. 

We now verify that the axioms of involutive $\infty$ sandwich cellularity hold for this choice of cell datum. The axiom (AC2) is satisfied since stacking diagrams can never increase the number of through lines, and the scalars which come from the removed closed components depend only on the top part of the diagram being stacked on. {Let $\hat{D}(\lambda,m)$ be the set of $(\lambda,m)$-diagrams with $\lambda$ through lines; it is potentially larger than $D(\lambda,m)$ as its diagrams are allowed to contain non-contractible loops.} We define $\Delta(\lambda)$ to be the free $k$-module with basis $\hat{D}(\lambda,m)$, which is naturally an $\overline{aPa_m}(\beta)$-$H_\lambda$-bimodule with action coming from the multiplication in $\overline{aPa_m}(\beta)$, {where we set $x\cdot y=0$ for $x\in \overline{aPa_m}(\beta),y\in \Delta(\lambda)$ if $xy$ has fewer than $\lambda$ through lines.} {By construction, every basis diagram in $\Delta(\lambda)$ factors uniquely as $t\cdot m$, where $t\in \mathcal{T}(\lambda), m\in B_\lambda$, so it is a free right $H_\lambda$-module.} Similarly, we define $\nabla(\lambda)$ to be the $H_\lambda$-$\overline{aPa_m}(\beta)$-bimodule analogously defined, which as a free $k$-module has basis $(\hat{D}(\lambda,m))^\star$ {and is a free left $H_\lambda$-module}.  Then $\Delta(\lambda)\otimes_{H_\lambda}\nabla (\lambda)$ is evidently isomorphic to $A^{\ge\lambda}/A^{>\lambda}$ as $\overline{aPa_m}(\beta)$-$\overline{aPa_m}(\beta)$-bimodule. Finally, the axiom (AC4) is clear from the diagrammatic definition of $(\rule{0.5em}{0.4pt}
)^\star$.

The above argument works for all the algebras which we have defined, so these are all involutive {$\infty$} sandwich cellular. Note that in the case where the diagrams involved are Temperley--Lieb or Brauer, the perfect matching condition restricts the possible number of through lines, and hence $\Lambda=\Lambda_2$ in these cases. We now discuss the algebras $H_\lambda$. For the planar affine algebras and $\lambda > 0$, non-contractible loops cannot appear (cf. \autoref{planar no loop}), so $H_\lambda$ is exactly $k[t,t^{-1}]\cong k[\Z]$ with the identification as in \autoref{Z idenfification}; when $\lambda=0$, there are no through lines but an arbitrary (finite) number of non-contractible loops is allowed, so we get $H_\lambda\cong k[t]$. For the reduced affine planar algebras, the only change is that  we have $H_0\cong k$ since the non-contractible loops are removed. For the planar periodic algebras, $H_\lambda=k$ if $\lambda>0$, since the only $(\lambda,\lambda)$-diagram with $\lambda$ through lines is identity, and $H_{0}=k[t]$ as before. For the reduced planar periodic diagrams we always have $H_{\lambda}=k$. For the symmetric affine diagram algebras, we have {$H_{\lambda}\cong k[\Z \wr S_{\lambda}]\otimes k[t]\cong k[\Z \wr S_{\lambda}][t]$} since the $(\lambda,\lambda)$-diagrams with $\lambda$ through lines generate the group algebra $k[\Z \wr S_\lambda]$ and we now allow non-contractible loops in diagrams with any number of through lines. Similarly, the periodic symmetric affine diagram algebras we get {$H_{\lambda}=k[\widetilde{S}_\lambda][t]$}. For the reduced versions, {the variable $t$ disappears}.  

With respect to these $\infty$ sandwich cellular data, the algebras associated with planar diagram algebras and the rook algebra are {sandwich} cellular, because the top sets involved are finite. The other algebras do not have finite top sets, since they have (in particular) cups which are allowed to have arbitrarily large winding number.

\renewcommand{\arraystretch}{1.5}
\begin{table}[H]
\centering
\resizebox{0.92\textwidth}{!}{
\begin{tabular}{l|c|c|c || l|c|c|c}
Algebra & $\Lambda$ & $H_\lambda$ & Sandwich &
Algebra & $\Lambda$ & $H_\lambda$ & Sandwich \\
\hline\hline

$\overline{aTL_m}$ & $\Lambda_2$ &
\multirow{2}{*}{$k[t,t^{-1}]$ if $\lambda>0$, $k[t]$ if $\lambda=0$} &
\multirow{8}{*}{yes} &

$\overline{aBr}_m$ & $\Lambda_2$ &
\multirow{2}{*}{{$k[\mathbb{Z}\wr S_\lambda][t]$}} &
\multirow{8}{*}{no} \\
\cline{1-2}\cline{5-6}

$\overline{aMo_m}$ & $\Lambda_1$ & & &
$\overline{aRoBr_m}/\overline{aPa_m}$ & $\Lambda_1$ & & \\

\cline{1-3}\cline{5-7}

$aTL_m$ & $\Lambda_2$ &
\multirow{2}{*}{$k[t,t^{-1}]$ if $\lambda>0$, $k$ if $\lambda=0$} & &
$aBr_m$ & $\Lambda_2$ &
\multirow{2}{*}{$k[\mathbb{Z}\wr S_\lambda]$} & \\

\cline{1-2}\cline{5-6}

$aPRo_m/aMo_m$ & $\Lambda_1$ & & &
{$aRo_m/aRoBr_m/aPa_m$} & $\Lambda_1$ & & \\

\cline{1-3}\cline{5-7}

$pPRo_m/pMo_m$ & $\Lambda_1$ &
\multirow{2}{*}{$k$} & &
$\overline{pBr_m}$ & $\Lambda_2$ &
\multirow{2}{*}{{$k[\widetilde{S}_\lambda][t]$}} & \\

\cline{1-2}\cline{5-6}

$pTL_m$ & $\Lambda_2$ & & &
$\overline{pRoBr_m}/\overline{pPa_m}$ & $\Lambda_1$ & & \\

\cline{1-3}\cline{5-7}

$\overline{pTL}_m$ & $\Lambda_2$ &
\multirow{2}{*}{$k$ if $\lambda>0$, $k[t]$ if $\lambda=0$} & &
$pBr_m$ & $\Lambda_2$ &
\multirow{2}{*}{$k[\widetilde{S}_\lambda]$} & \\

\cline{1-2}\cline{5-6}

$\overline{pMo}_m$ & $\Lambda_1$ & & &
$pRo_m/pRoBr_m/pPa_m$ & $\Lambda_1$ & & \\

\cline{1-3}\cline{5-7}

\hline
\end{tabular}}

\caption{($\infty$) sandwich cellular data for the affine and periodic algebras, arranged in two columns. `Sandwich' means honest sandwich and not $\infty$ sandwich.}
\label{tab:sandwich}
\end{table}

\begin{Theorem}
Over a field $k$, the diagram algebras in \autoref{tab:sandwich} are involutive ($\infty$) sandwich cellular with the specified data.
\end{Theorem}

\begin{proof}
This is discussed above.
\end{proof}

\subsection{Simple modules}

We now use the theory of sandwich cellular algebras to construct all simple modules for the algebras which are involutive sandwich cellular. For background we refer to \cite[\S 2]{tubbenhauer2022handlebody}. Recall that for an involutive sandwich cellular {algebra} $A=(\Lambda, \mathcal{T},(H_\lambda,B_\lambda, {(\rule{0.5em}{0.4pt})^{\star_\lambda}}),C, (\rule{0.5em}{0.4pt} 
)^\star)$, the \textbf{apex} of an $A$-module, if it exists, is a maximal $\lambda \in \Lambda$ such that the submodule $k\{c_{T,m,B}^\lambda: T,B \in \mathcal{T}(\lambda), m\in H_\lambda\}$ does not annihilate $M$. Every simple module has an apex. For $\lambda \in \Lambda$, the multiplication on (either side of) $\Delta(\lambda)\otimes_{H_\lambda}\nabla(\lambda)$ is determined by a bilinear form valued in the algebra $H_{\lambda}$, \[\phi^\lambda:\nabla(\lambda)\otimes_{A}\Delta(\lambda)\to H_{\lambda},\]
which in our case can be thought of as follows: given $x\in \Delta(\lambda)$ and $y^\star \in \nabla(\lambda)$, the pairing $\phi^\lambda(y^\star \otimes x)$ is just $y^\star x$ as an element of $H_\lambda$ ({if it has $\lambda$ through lines; otherwise it is 0}). This extends to a bilinear form {\[\phi^\lambda_S=\phi^\lambda \otimes_{H_{\lambda}} id_S:  \nabla(\lambda)\otimes_{A}\Delta(\lambda)\otimes_{H_{\lambda}} S \to H_{\lambda}\otimes_{H_\lambda} S\cong S, \]}
where $S$ is any simple module of $H_{\lambda}$. The bilinear form $\phi_S^\lambda$ induces an associated map {$\overline{\phi}_S^\lambda: \Delta(\lambda)\otimes_{H_\lambda} S\to \Hom_{H_\lambda}(\nabla(\lambda),H_{\lambda})\otimes_{H_{\lambda}} S$}, and we write $\rad(\lambda, S)$ for the kernel of the map. We also write $\Delta(\lambda, S)$ for the $A$-module $\Delta(\lambda)\otimes_{H_\lambda} S$. The quotient $\Delta(\lambda, S)/\rad(\lambda, S)$ is a simple $A$-module, which we denote by $L(\lambda, S)$. Below we write $\Lambda$ to mean either $\Lambda_1$ or $\Lambda_2$ depending on the algebra, cf. \autoref{tab:sandwich}. If $G$ is a group, we write $Irr_k(G)$ for the set of simple $kG$-modules up to equivalence.

\begin{Theorem}\label{irrep classification theorem}
Let $k$ be an algebraically closed field.
\begin{enumerate}
\item The simple modules of $\overline{aTL_m}(\beta)$ and $\overline{aMo_m}(\beta)$ are parameterized by $\{(\lambda,z): \lambda\in \Lambda, z \in k^\times\}$ with an additional simple module indexed by $(0,0)$ if $0\in \Lambda$ {and $\beta\neq 0$}. {If $\beta=0$, then for $\overline{aMo_m}(0)$ all parameters $(\lambda,z)$ with $m-\lambda$ odd are removed.}

\item The simple modules of $aTL_m(\beta,\alpha)$, $aMo_m(\beta,\alpha)$ and $aPRo_m(\beta)$ are parameterized by $\{(\lambda,z):\lambda\in\Lambda\setminus\{0\}, z\in k^\times\}$ with an additional simple module indexed by $(0,1)$ if $0\in\Lambda$. {If $\beta=0$, then for $aMo_m(0,\alpha)$ all parameters $(\lambda,z)$ with $m-\lambda$ odd are removed, while for $aPRo_m(0)$ all parameters $(\lambda,z)$ with $\lambda<m$ are removed. Moreover, if $\beta=\alpha=0$, the parameter $(0,1)$ is removed for $aTL_m(0,0)$ and $aMo_m(0,0)$.}

\item The simple modules of $\overline{pTL_m}(\beta)$ and $\overline{pMo_m}(\beta)$ are parameterized by $\{\lambda:\lambda\in\Lambda\setminus\{0\}\}$ with additional simple modules indexed by $\{(0,z):z\in k\}$ if $0\in\Lambda$, {where $(0,0)$ is removed if $\beta=0$}. {If $\beta=0$, then for $\overline{pMo_m}(0)$ all parameters $\lambda\in\Lambda\setminus\{0\}$ with $m-\lambda$ odd are removed. }

\item The simple modules of $pTL_m(\beta,\alpha)$, $pMo_m(\beta,\alpha)$ and $pPRo_m(\beta)$ are parameterized by $\{\lambda:\lambda\in\Lambda\}$. {If $\beta=0$, then for $pMo_m(0,\alpha)$ all parameters $\lambda$ with $m-\lambda$ odd are removed, while for $pPRo_m(0)$ all parameters $\lambda$ with $\lambda<m$ are removed. Moreover, if $\beta=\alpha=0$, the parameter $0$ is removed for $pTL_m(0,0)$ and $pMo_m(0,0)$.}

\item The simple modules of $aRo_m(\beta)$ are parameterized by $\{(\lambda,S): \lambda\in \Lambda, S \in Irr_k(\Z \wr S_\lambda)\}$. Over $\C$, $Irr_k(\Z \wr S_\lambda)$ is $\Gamma_\lambda$ defined in \autoref{rep of wreath}. {If $\beta=0$, then all parameters $(\lambda,S)$ with $\lambda<m$ are removed.}

\item The simple modules of $pRo_m(\beta)$ are parameterized by $\{(\lambda,S):\lambda\in \Lambda, S \in Irr_k(\widetilde{S}_\lambda)\}$. {If $\beta=0$, then all parameters $(\lambda,S)$ with $\lambda<m$ are removed.}
\end{enumerate}
Moreover, the simple modules are all of the form
\[
\Delta(\lambda,S)/\rad(\lambda,S)=L(\lambda,S),
\]
where {$(\lambda,S)$ ranges over the parameters listed above.}
\end{Theorem}
\begin{proof}
By \cite[Theorem 2.16 (c)]{tubbenhauer2022handlebody}, {if $\lambda\in \Lambda$ is an apex,} there is a bijection \[\{\text{Simple\ } A\text{-modules with apex\ }\lambda \}\leftrightarrow \{\text{Simple\ } H_\lambda\text{-modules}\}\] given by $L(\lambda, S)\leftrightarrow S$. ({To be precise, this requires the technical hypothesis of \cite[Lemma 2.11]{tubbenhauer2022handlebody}\footnote{This lemma contains a typo in the statement and has an unstated hypothesis that the element $F\in H_\lambda$ in its statement should be invertible, which is used in its proof. When $\beta=\lambda=0$ one cannot choose $F$ to be invertible for the affine Temperley--Lieb algebra: diagrammatically this corresponds to the fact that non-contractible loops cannot be removed.} which fails for the Temperley--Lieb and Motzkin algebras when $\lambda=\beta=0$, and for their reduced versions when $\lambda=\beta=\alpha=0$. In this case it is still true, by \cite[Lemma 2.15(c)]{tubbenhauer2022handlebody} that the simples with apex $0$ correspond to the simples of the algebra $A_0$ defined in \eqref{Alambda}; in the present case they can be analyzed using \cite[Theorem 3.10]{koenig2012affine}. For example, in the affine Temperley--Lieb case, when $\beta=0$ the simples with apex 0 are parameterized by $k^\times$. This is the source of the exceptions.)} By \cite[Theorem 2.16 (a)]{tubbenhauer2022handlebody}, an element $\lambda \in \Lambda$ is an apex if and only if $\phi^\lambda \not\equiv 0$. {An easy computation, cf. \cite[Remarks (2.7)]{graham1998representation}, shows that for $\overline{aTL_m}(\beta)$ and $\overline{pTL_m}(\beta)$ the form never vanishes, and for their reduced versions the form only vanishes if $\lambda=\beta=\alpha=0$}. {If $\beta\neq 0$,} for $aRo_m(\beta),pRo_m(\beta), aPRo_m(\beta)$ and $pPRo_m(\beta)$ the form $\phi^\lambda$ is never 0 (in fact, the restriction of $\phi^\lambda$ to the submodule spanned by ordinary planar rook diagrams always has its Gram matrix {equal to $\beta^{m-\lambda}$ times} a permutation matrix, see e.g. \cite[Proposition 4F.7]{KhSiTu-monoidal-cryptography}). {When $\beta=0$, $\phi^\lambda$ vanishes whenever $\lambda<m$, since to preserve the number of through lines any isolated vertex must be paired with another one, which creates a factor of $\beta=0$.} For the algebras associated to Motzkin diagrams {the analysis is by combining the two preceding cases: when $m-\lambda$ is odd, there must be isolated vertices (as opposed to only cups) in all top diagrams, so the form vanishes if and only if $\beta=0$. When $m-\lambda$ is even, there are top diagrams which are Temperley--Lieb diagrams, and hence the form does not vanish unless $\lambda=\alpha=\beta=0$}. The claims now follow by considering the simple modules of the algebras $H_\lambda$ case by case.  
\end{proof}

\begin{Remark} \label{GL parameterization}
In \cite{graham1998representation}, the simple modules for $\overline{aTL}_m(\beta)$ are indexed by the set $\{(\lambda, z):\lambda \in \Lambda_2, z\in k^\times\}/\sim$, where $\sim$ is the equivalence relation identifying $(0,z)\sim (0,z^{-1})$, with $(0,q)\sim (0,q^{-1})$ removed if $\beta=0,m\in 2\Z$ and $q+q^{-1}=0$. This is equivalent to our parameterization because {$k^\times/(z\sim z^{-1})\to k, z\mapsto z+z^{-1}$} is a bijection when $k$ is algebraically closed.
\end{Remark}

\begin{Proposition}\label{planar rook reducible}
{Over any field $k$ and for $\beta\neq 0$, the cell modules $\Delta(\lambda, S)$ of $aPRo_m(\beta)$ are simple, i.e. $\Delta(\lambda, S)=L(\lambda, S)$.}   
\end{Proposition}
\begin{proof}
For $aPRo_m(\beta)$ the top sets and diagram stacking rules are the same as those for the ordinary planar rook algebra $PRo_m(\beta)$, {so they also have the same Gram matrices for the bilinear forms valued in $H_\lambda$. In the planar rook case this matrix is known to be given by some power of $\beta$ multiplied by identity, see \cite[Proposition 4F.7]{KhSiTu-monoidal-cryptography}}. {It follows that after extension to take value in $S$, the bilinear form is still non-degenerate, which implies the result.}    
\end{proof}

\subsection{Cardinalities of top sets}
In this section we determine the cardinalities of the top sets $\mathcal{T}(\lambda)$ for the algebras that have finite top sets. {Assuming we work over an algebraically closed field,} for the planar algebras $|\mathcal{T}(\lambda)|$ is the same as the dimension of the cell modules. For affine or periodic rook algebra, the dimension of $\Delta(\lambda, S)$ is obtained by multiplying $\lvert \mathcal{T}(\lambda)\rvert$ with $\dim S$. We use square brackets to mean `coefficient of'.

\begin{Lemma}[Lagrange inversion]\label{lagrange lemma}
Let $F(x)$ be a formal power series with $F(0)=0$, $F(x)=x\phi(F(x))$, where $\phi(t)$ is a formal power series with $\phi(0)\neq 0$. Then for any formal power series $G(t)$ and $n\ge 1$, \[[x^n]G(F(x))=\frac{1}{n}[t^{n-1}]G'(t)(\phi(t))^n.\]
{Here $G'$ denotes the formal derivative of $G$.}
In particular, for $G(t)=t^k$, $n\ge k\ge 1$, \[[x^n](F(x))^k=\frac{k}{n}[t^{n-k}](\phi(t))^n.\]
\end{Lemma}
\begin{proof}
See e.g. \cite[Theorem 2.1.1]{GESSEL2016212}.    
\end{proof}

\begin{Proposition}
The cardinalities of the top sets $\mathcal{T}(\lambda)$ are:
\begin{enumerate}
\item For $\overline{aTL_m}(\beta)$, $\overline{pTL_m}(\beta)$ or their reduced versions: $\binom{m}{(m-\lambda)/2}$.
\item For $\overline{aMo_m}(\beta)$, $\overline{pMo_m}(\beta)$ or their reduced versions: \[ \sum_{j=0}^{\floor{(m-\lambda)/2}{}} \frac{m!}{j!(m-\lambda-2j)!(\lambda+j)!}=\sum_{j=0}^{\floor{(m-\lambda)/2}{}} \binom{m}{j}\binom{m-j}{m-\lambda-2j} =[t^{m-\lambda}](1+t+t^2)^m.\]
\item For {$aPRo_m(\beta)$ or $pPRo_m(\beta):\binom{m}{\lambda}$}. 
\item For $aRo_m(\beta)$, $pRo_m(\beta)$ or their reduced versions: $\binom{m}{\lambda}$. 
\end{enumerate}
\end{Proposition}

\begin{proof}
  {For the top sets it does not matter whether we consider the affine or periodic versions.} The affine Temperley--Lieb case was proven in \cite[Corollary 1.12]{graham1998representation} by induction; alternatively a counting proof can be given by adapting the argument below. In the Motzkin case, to give an element of $\mathcal{T}(\lambda)$ we first specify the positions of $\lambda$ through lines, which can be thought of as dividing a circle into $\lambda$ segments with lengths $l_1,\dots, l_\lambda$, $\sum l_i=m-\lambda$. {Assume first $\lambda>0$.} The number of ways to fill these segments with non-intersecting cups and isolated components is $\prod_i M_{l_i}$ where $M_{l_i}$ is the $l_i$-th Motzkin number. We get \[\lvert\mathcal{T}(\lambda) \rvert=\frac{m}{\lambda}\sum_{l_1+\dots+l_\lambda=m-\lambda}\prod_{i=1}^\lambda M_{l_i}=\frac{m}{\lambda}[x^{m-\lambda}]M(x)^\lambda,\]
where $M(x)=1+x+2x^2+4x^3+\dots$ is the generating function of Motzkin numbers satisfying the relation $M(x)=1+xM(x)+x^2M(x)^2$. Write $\phi(t)=1+t+t^2$, $N(x)=xM(x)$, then using \autoref{lagrange lemma} we get $$\frac{m}{\lambda}[x^{m-\lambda}]M(x)^\lambda = \frac{m}{\lambda}[x^m]N(x)^\lambda=[t^{m-\lambda}](1+t+t^2)^m,$$
which expands to give the desired formula. {In the $\lambda=0$ case, note that Motzkin diagrams on a circle are in bijection with words in $\{-1,0,1\}^m$ having total sum 0, and thus $|\mathcal{T}(0)|=[x^0](x^{-1}+1+x)^m=[x^m](1+x+x^2)^m$, since $(x^{-1}+1+x)^m=x^{-m}(1+x+x^2)^m$.} In the planar rook case the top sets are the same as in the ordinary case, so the formula follows from \cite[\S 2]{flath2009planar}. Finally, a top rook diagram is just a top planar rook diagram. 
The proof is complete.
\end{proof}

\subsection{Finite dimensional quotients}\label{fin dim section}
The algebras ${aPa_m}(\beta),{aBr_m}(\beta), {aRoBr_m}(\beta)$ are not sandwich cellular, because the `top sets' are infinite (the cups are allowed to wind around an arbitrary amount). However, we can consider their finite-dimensional quotients, which are sandwich cellular.
\begin{Definition}
For $r\in \Z_{\ge 1}$, denote by $aPa_{m,r}(\beta)$ (resp. $aBr_{m,r}(\beta)$, $aRoBr_{m,r}(\beta)$) the $r$-reduced affine partition (resp. Brauer, rook-Brauer) algebra, obtained from $aPa_m(\beta)$ (resp.$aBr_m(\beta)$, $aRoBr_m(\beta)$) by imposing the additional relations $t_i^r=1,1\le i \le m$.  
\end{Definition}
Note that when $r=1$, we recover the non-affine diagram algebras.

\begin{Proposition}
The algebras $aPa_{m,r}(\beta)$, $aBr_{m,r}(\beta)$ and $aRoBr_{m,r}(\beta)$ are involutive sandwich cellular algebras, with the same poset $\Lambda$ and anti-isomorphism $(\rule{0.5em}{0.4pt}
)^\star$ as before, {the top sets} $\mathcal{T}_{m,r}(\lambda)$ given {similarly by picking orbit representatives}, and $H_\lambda$ isomorphic to the group algebra {$kS(r,\lambda)$}, where {$S(r,\lambda)$} is the generalized symmetric group {$C_r\wr S_\lambda$}, with a basis $B_\lambda$ given by the invertible diagrams.
\end{Proposition}

\begin{proof}
Similar to the proof for the infinite-dimensional algebras. 
\end{proof}

\begin{Proposition}\label{r affine irrep} Let $k$ be a field.
The simple modules of $kaPa_{m,r}(\beta)$, $kaRoBr_{m,r}(\beta)$, and $kaBr_{m,r}(\beta)$ over $k$ are parameterized by $\{(\lambda, S): \lambda \in \Lambda, S\in {Irr}_k(S(r,\lambda)), {\phi^\lambda\not\equiv 0}\}$ where $\Lambda=\Lambda_1$ for the first two algebras, and $\Lambda=\Lambda_2$ for the last one.
\end{Proposition}
\begin{proof}
Similar to that of \autoref{irrep classification theorem}. 
\end{proof}

\begin{Remark}
By inflating the simple modules described in \autoref{r affine irrep} for $r\ge 1$, we obtain all simple modules of the symmetric reduced affine algebras on which $t_i$ acts as an operator of finite order.   
\end{Remark}

\begin{Proposition}\label{top sets finite}
The cardinalities of the top sets $\mathcal{T}_{m,r}(\lambda)$ are:
\begin{enumerate}
\item For $aBr_{m,r}(\beta)$: \[\lvert \mathcal{T}_{m,r}(\lambda)\rvert =\binom{m}{\lambda}{r}^{(m-\lambda)/2}(m-\lambda-1)!!.\]
\item For $aRoBr_{m,r}(\beta)$:\[\lvert \mathcal{T}_{m,r}(\lambda)\rvert= \binom{m}{\lambda}\sum_{p=0}^{\floor{(m-\lambda)/2}{}}\binom{m-\lambda}{2p}(2p-1)!!{r}^p.\]
\item For $aPa_{m,r}(\beta)$:
\[\lvert \mathcal{T}_{m,r}(\lambda)\rvert ={\sum_{t=\lambda}^m\binom{t}{\lambda}\sigma(m,t)r^{m-t}},\]
where $\sigma(m,t)$ is the Stirling number of the second kind.
\end{enumerate}
\end{Proposition}
\begin{proof} For $aBr_{m,r}(\beta)$, we first choose the positions of $\lambda$ {top vertices connected to} through lines. For the remaining $m-\lambda$ vertices in the fundamental rectangle, there are $(m-\lambda-1)!!$ ways to pair them in the fundamental rectangle, and for each pair a choice of displacement {$\delta \in {\{0,\dots, r-1\}}$, so $r^{(m-\lambda)/2}$ choices altogether. This gives the formula.} The case for $aRoBr_{m,r}(\beta)$ is similar: in this case we first choose $2p$ of the $m-\lambda$ `non-through' vertices to pair with each other, so that the remaining vertices are isolated. For $aPa_{m,r}(\beta)$, {we first partition the $m$ top vertices in the fundamental rectangle into $t$ blocks; $\lambda$ of these blocks are connected to a bottom vertex. Ignoring displacements for the moment, there are $\sigma(m,t)$ such set partitions, where $\sigma(m,t)$ is the Stirling number of the second kind. If one of these blocks has $q$ vertices, then it has $r^q/r=r^{q-1}$ possible configurations modulo $r$ (note that it does not matter whether the block is one of the $\lambda$ `through' blocks and hence has a bottom vertex, since we are taking orbit representatives). Now suppose the $t$ blocks have $q_1,\dots, q_t$ vertices respectively, then the total number of configurations is $r^{q_1-1}\dots r^{q_t-1}=r^{m-t}$. This gives the result.}  
\end{proof}

Using \autoref{top sets finite}, we can deduce the dimensions of the cell modules using {$\dim \Delta(\lambda, S)=\lvert \mathcal{T}_{m,r}(\lambda)\rvert \cdot \dim S$.}

\section{Growth Problems}\label{growth section}

In this section we consider the underlying monoids of the algebras which we have defined, and study the tensor powers of their representations.
\subsection{General theory}

Throughout let $M$ be a potentially infinite monoid with group of units $G$, and let $k$ be an algebraically closed field {of characteristic $p\ge 0$}. 
{All $kM$-modules considered below are assumed to be finite-dimensional.} {Throughout this section, we also assume that $M$ is \textit{Dedekind-finite}, which means every left-invertible element is right-invertible; this is equivalent to $M\setminus G$ being a two-sided ideal. All finite monoids are Dedekind-finite, and so are all the infinite diagram monoids considered in this paper (because the number of through lines can never increase, and non-contractible loops cannot be removed.) In almost all results below we assume the existence of a $kM$-module $V$ on which $M\setminus G$ acts non-invertibly, which already implies Dedekind-finite.}

If $V$ is a $kM$-module, the corresponding \textit{fusion graph} for summand $\Gamma^b$ (resp. fusion graph for length $\Gamma^l$) is the (oriented and weighted) graph whose vertices are indecomposable (resp. simple) $kM$-modules which are summands (resp. composition factors) of $V^{\otimes n}$, for some $n\ge 0$, and {such} that there is an edge of weight $m$ from the vertex $V_i$ to the vertex $V_j$ if $V_j$ occurs $m$ times in the direct sum decomposition (resp. composition series) of $V\otimes V_i$. The corresponding (potentially countably infinite) adjacency matrix is called the \textit{action matrix} for summand (resp. for length). Let $\lambda_b^{\mathrm{sec}}$ (resp. $\lambda_l^{\mathrm{sec}}$) denote any second largest eigenvalue (in terms of modulus) of the corresponding action matrix, the largest always being the dimension {of $V$}. This is not always defined, but will be in the settings we are interested in, see \cite[Lemma 5.2]{lacabanne2024asymptoticsinfinitemonoidalcategories}. Let $\Res_G(V)$ denote the $kG$-module that comes from restricting the action on $V$; we write $Z_V(G)$ for the subgroup of elements in $G$ which act as scalars on $V$, and {for $g\in Z_V(G)$ we let $\omega_V(g)$ denote the scalar by which $g$ acts on $V$}.
Suppose that $G$ has $N$ ($p$-regular) conjugacy classes, with representatives $g_1,\dots, g_N$. Finally, we say $M$ satisfies the \textit{group-injective condition} (over $k$) if some {nonzero} projective $kG$ module $P$, made into a $kM$-module by letting $M\setminus G$ act as 0 (we write $\Ind_G(P)$ for this module), is injective as a $kM$-module; {note that $\Ind_G(P)$ only makes sense under the Dedekind-finite assumption.}

We write $f(n)\sim g(n)$ to mean asymptotically equal, i.e. $f(n)/g(n)\to 1$ as $n\to \infty$. We also use the usual capital $\mathcal{O}$ notation. The following is known about the growth problems for finite monoids:

\begin{Theorem} \label{finite theorem} Let $M$ be a finite monoid. Suppose $V$ is a $kM$-module on which no non-identity element acts {as identity}. If $M$ satisfies the group-injective condition, then
\begin{enumerate}
 \item We have \begin{equation} b(n)\sim a(n):=\frac{1}{\lvert G\rvert}\sum_{\substack{1\le t\le N \\ g_t \in Z_V(G)}} S_{t} {\cdot}  \omega_V(g_t)^n\cdot (\dim V)^n, 
    \end{equation}
    where $S_t$ is the sum over entries of the column corresponding to $g_t^{-1}$ in the (Brauer) character table of irreducible characters. In other words, $V$ has the same asymptotic growth rate as the $kG$-module $\Res_G(V)$, cf. \cite[Theorem 1]{he2025growthproblemsrepresentationsfinite}.
\item $\lvert b(n)/a(n)-1\rvert \in \mathcal{O}(\lvert \lambda_b^{\mathrm{sec}}/\dim V \rvert^n + n^{-c})$, for some constant $c>0$, and
    \item $\lvert b(n)-a(n) \rvert \in \mathcal{O}(\lvert \lambda_b^{\mathrm{sec}}\rvert^n + n^{d})$ for some constant $d > 0$.
\end{enumerate}
Exactly analogous results hold for $l(n)$ (without needing the assumption that $M$ satisfies the group-injective condition). In particular, we have
\[l(n)\sim k(n):=\frac{1}{\lvert G\rvert}\sum_{\substack{1\le t\le N \\ g_t \in Z_V(G)}} T_{t}{\cdot} \omega_V(g_t)^n\cdot (\dim V)^n,\]
where $T_t$ is the sum over entries of the column corresponding to $g_t^{-1}$ in the (projective) Brauer character table.
\end{Theorem}
\begin{proof}
See \cite[Theorem 1]{he2025growthproblemsrepresentationsfinite2} and \cite[Theorem 2B.2]{he2025tensorpowersrepresentationsdiagram} for the results for $b(n)$ and $l(n)$ respectively.
\end{proof}

When $M$ is a finite monoid (or semigroup), we always have $C\cdot (\dim V)^n \le b(n)\le l(n)\le (\dim V)^n$ for some constant $C$, see \cite[Proposition 1]{he2025growthproblemsrepresentationsfinite2}. When $M$ is infinite, this need not be true, as illustrated in the following example.

\begin{Example}\label{jordan example}
An indecomposable representation of $\Z$ (or $\N$) is a Jordan block. Let $J_m(1)$ denote the unipotent Jordan block of size $m$. By \cite{marcus1975elementary} we have that over $\C$, $$J_m(1)\otimes J_p(1)=\bigoplus_{i=1}^{\min(m,p)}J_{m+p-2i+1}(1), $$
which is exactly the Clebsch-Gordan formula for simple modules of $SL(2,\C)$. {Applying} \cite[Example 2.8]{coulembier2023asymptotic}, we get 
$$b(n)\sim \sqrt{\frac{6}{\pi(m^2-1)}}\frac{m^n}{\sqrt{n}},$$
{which has $n^{-1/2}$ as a correction factor instead of plain exponential growth.}
\end{Example}

However, under certain assumptions the behavior of $b(n)$ and $l(n)$ for infinite monoids is exactly like in the finite case. 

\begin{Definition} \label{monoid def} Let $M$ be a monoid with group of units $G$.
\begin{enumerate}
\item We say $M$ is \textit{right-finite} if there is a finite set $X$ such that any $m\in M\setminus G$ can be written as a word ending in some $x\in X$. In particular, finite monoids are right-finite.
\item We say $M$ is \textit{right idempotent generated} if every $m\in M\setminus G$ can be written as $m=ne$, where $e\neq 1$ is an idempotent. We say $M$ is (von Neumann) \textit{regular} if every $m\in M$ can be written as $m=mam$ for some $a\in M$. If $M$ is regular then we can always write $m=m(am)$, where $am$ is an idempotent, so $M$ is right idempotent generated.
\item We say $M$ is \textit{dimension-bounded} (over $k$) if all simple $kM$-modules are bounded in dimension by some global constant. In particular, finite monoids are dimension-bounded over any field.
\end{enumerate}
We also extend these definitions to the case where $M$ is a semigroup $S$ without unit, by setting $G=\emptyset$.
\end{Definition}

{Note that right-finite does not imply Dedekind-finite, as the bicyclic monoid $M=\langle p,q\mid pq=1\rangle$ shows.}

\begin{Lemma}\label{generalised length}
Let $M$ be right-finite, and let $V$ be a $kM$-module on which no nonunit acts invertibly. Then $V^{\otimes (\lvert X\rvert)}$ has a nonzero submodule on which $M\setminus G$ acts as 0, where $X$ is the finite set in the definition of right-finiteness. 
\end{Lemma}
\begin{proof}
For each $x\in X$, pick $v_x$ in the null space of the action of $x$, then $\bigotimes v_x \in V^{\otimes (\lvert X\rvert)}$ is a nonzero vector in the submodule of elements on which nonunits act as 0. {Note that if $W$ is any $kM$-module, the set of elements $w\in W$ on which $M\setminus G$ acts as zero is always a submodule, as can be easily seen.}    
\end{proof}

\begin{Remark}
If all nonunits $m$ of $M$ {generate a finite cyclic subsemigroup $\{m,m^2,\dots\}$ (e.g. if $M$ is finite)}, requiring that no nonunit acts invertibly is the same as requiring that no nonunit acts as identity, and is thus a strictly weaker condition than faithfulness. If $M$ is an infinite monoid, faithful modules need not satisfy the condition: for example, a non-nilpotent Jordan block is a faithful module of $\N$ with all nonunits acting invertibly. Note that \autoref{jordan example} shows that the condition that no nonunit acts invertibly is necessary. 
\end{Remark}

\begin{Definition} Let {$A$} be an action matrix (for summand or length) with fusion graph $\Gamma$. {Write $a_{ij}^n$ for the $(i,j)$-th entry of $A^n$.}
\begin{enumerate}
\item {If $A$ is irreducible (as a matrix),} the \textit{period} of {$A$} (or of $\Gamma$) is the greatest common divisor of all $n$ such that $a_{ii}^n \neq 0$, where $i$ {runs over the indexing set for the rows and columns}. 
\item {If $A$ is irreducible,} the (\textit{Vere--Jones) {\textup{PF} dimension}} of {$A$} (or $\Gamma$)
  is $\textup{PF}\dim^{vj} {A} := \lim\limits_{n\to \infty} \sqrt[h \cdot n]{a_{ij}^{h\cdot n}}$, where $h \in \Z_{\ge 1}$ is the period of $A$  {and the pair $(i,j)$ is such that the terms in the sequence are not identically zero. By \cite[Theorem A]{vere1967ergodic}, the limit exists and is independent of the choice of such $(i,j)$.} 
  \item A \textit{class} $C(\Gamma)$ of $\Gamma$ is a strongly connected component of $\Gamma$; it is called a \textit{basic class} if moreover $\textup{PF}\dim^{vj} C(\Gamma) \ge \textup{PF}\dim^{vj} C'(\Gamma)$ for all classes $C'(\Gamma)$, and is called a \textit{final basic class (FBC)} if in addition, there is no path from it to any other basic class. It is called \textit{final in $\Gamma$} if there is no path to any vertex outside $C(\Gamma)$.  
\end{enumerate}
(We refer to \cite[\S 4]{lacabanne2024asymptoticsinfinitemonoidalcategories} for more explanations of the definitions.)
\end{Definition}

\begin{Notation}
We denote by $\Gamma^l_G$ (resp. $\Gamma^b_G$) the subgraph of $\Gamma^l$ (resp. $\Gamma^b$) consisting of all simple (resp. projective indecomposable) $kM$-modules on which $M\setminus G$ acts as 0, i.e. modules of the form $\Ind_G(V)$ where $V$ is a {simple (resp. projective indecomposable)} $kG$-module.
\end{Notation}

It is conjectured in \cite[Conjecture 1]{he2025growthproblemsrepresentationsfinite2} that \autoref{finite theorem} is true for a finite monoid without the hypothesis that $M$ satisfies the group-injective condition over $k$. In particular, in \cite{he2025tensorpowersrepresentationsdiagram} the conjecture is verified to be true over $\C$ for all the diagram monoids in \autoref{fig: table of algebras}. We will now work towards proving a generalization of \autoref{finite theorem} which in particular will imply \cite[Conjecture 1]{he2025growthproblemsrepresentationsfinite2} under the hypothesis that $M$ is right idempotent generated.

\begin{Lemma}\label{exponential factor lemma}
Let $V$ be a $kM$-module. Then 
\[\lim\limits_{n\to \infty}\sqrt[n]{b(n)}=\lim\limits_{n\to \infty}\sqrt[n]{l(n)}=\dim V.\]
\end{Lemma}
\begin{proof}
See \cite[Theorem 1.4]{Coulembier_2023}.
\end{proof}

\begin{Lemma}\label{eventually all group}
Suppose $M$ is right-finite and $V$ is a $kM$-module on which nonunits act non-invertibly. Suppose $M$ is moreover dimension-bounded, say by $N$. Write $r_n(V)$ for the fraction of composition factors of $V^{\otimes n}$ on which $M\setminus G$ acts as 0, then $r_n(V)$ converges to 1 exponentially as $n\to \infty$.
\end{Lemma}
\begin{proof}
We note that by \autoref{generalised length} any vertex $v$ in $\Gamma^l$ {is connected to a vertex in $\Gamma^l_G$ by a path of length at most $m:=|X|$, where $X$ is the finite set from \autoref{monoid def}}. {There can be no edges leaving $\Gamma_G^l$ (see proof of \autoref{inf theorem} below), so we can extend the path and assume it has length exactly $m$.} {Write $\Pi_n$ for the multiset of length $n$ paths starting from the trivial module (thinking of a weighted edge as multiple parallel edges), and write $A_n$ to be the subset of $\Pi_n$ of paths that avoid $\Gamma^l_G$. It suffices to show that $\lvert A_n \rvert / \lvert \Pi_n \rvert \to 0$ exponentially as $n\to \infty$.}

{For a path $p$, let $e(p)$ denote its end vertex. Since $1\le \dim e(p)\le N$, we get
$|\Pi_n|\le \sum_{p\in \Pi_n}\dim e(p)\le N|\Pi_n|$
and the same inequalities hold with $\Pi_n$ replaced by $A_n$. Write $d=\dim V$. The above implies 
\[\frac{|A_n|}{|\Pi_n|}\le \frac{\sum_{p\in A_n}\dim e(p)}{\frac{1}{N}\sum_{q\in \Pi_n}\dim e(q)}=\frac{N\sum_{p\in A_n}\dim e(p)}{d^n}.\]}
Put $s=n-rm$, where $r=\floor{n/m}{}$ {and $0\le s<m$ is the residue}. {For $0\le j<r$ and $p\in A_{jm}$, the sum $\sum \dim e(q')$ over $q'\in A_{(j+1)m}$ extending $p$ is bounded above by $\dim e(p)d^m-1\le (1-1/Nd^m)\dim e(p)d^m. $ Summing over the $p$'s and iterating gives 
\[\frac{\sum_{p\in A_{rm}}\dim e(p)}{d^{rm}}\le \left(1-\frac{1}{Nd^m}\right)^r.\]
Now we have 
\[\frac{|A_n|}{|\Pi_n|}\le \frac{N\sum_{p\in A_n}\dim e(p)}{d^n} \le  \frac{d^s\cdot N\sum_{p\in A_{rm}}\dim e(p)}{d^{rm}\cdot d^s}\le N\left(1-\frac{1}{Nd^m}\right)^{\floor{n/m}{}}\]
which tends to $0$ exponentially as $n\to \infty$, completing the proof.} 
\end{proof}

\begin{Lemma}\label{all comp factors group}
Let $M$ be right-finite and dimension-bounded, and let $V$ be a $kM$-module on which all nonunits act non-invertibly. Then for some $N\ge 1$, $V^{\otimes N}$ contains a summand all of whose composition factors {have $M\setminus G$ acting as 0}.
\end{Lemma}
\begin{proof}
 Suppose for a contradiction that the contrary is true. We know from \autoref{eventually all group} that as $n\to \infty$, the ratio $r_n(V)$ approaches 1, so since each indecomposable summand of $V^{\otimes n}$ contains a composition factor which {has some nonunit not acting as 0}, we have $b(n)\le l(n)(1-r_n(V))=l(n)-l^G(n)$, where $l^G(n)$ is the number of composition factors of $V^{\otimes n}$ which {have all nonunits acting as 0}. However, the inequality $b(n)/l(n)\le 1-r_n(V)$ leads to a contradiction since
 \[1=\limsup\limits_{n\to \infty} \sqrt[n]{b(n)/l(n)}>\limsup\limits_{n\to \infty} \sqrt[n]{1-r_n(V)}\]
  by \autoref{exponential factor lemma} and \autoref{eventually all group}. The proof is complete.    
\end{proof}

\begin{Theorem}\label{inf theorem}
Let $M$ be right-finite. Let $V$ be a $kM$-module on which nonunits act noninvertibly, and such that $\Res_G(V)$ is faithful. 
\begin{enumerate}
\item If $G$ is finite, then the results of \autoref{finite theorem} for $l(n)$ hold. If $G$ is finite and the group-injective condition is satisfied (or more generally, assume that for some $N$, $V^{\otimes N}$ contains a summand {with nonunits acting as 0}), then the results of \autoref{finite theorem} for $b(n)$ also hold. 
\item If $G$ is abelian, or more generally all simple $kG$-modules are one-dimensional, then the results of \autoref{finite theorem} for $l(n)$ hold. In particular, we have $l(n)\sim (\dim V)^n$.
\end{enumerate}
Suppose that $M$ is moreover dimension-bounded and right idempotent generated.  
\begin{enumerate}[resume]
 \item $V^{\otimes N}$ contains a summand {on which $M\setminus G$ acts as 0,} for some $N$.
    \item If $G$ is finite, then the results of \autoref{finite theorem} for $b(n)$ hold. In particular, this proves \cite[Conjecture 1]{he2025growthproblemsrepresentationsfinite2} when $M$ is right idempotent generated.  
\end{enumerate}
  
\end{Theorem}

\begin{proof}
(a). \autoref{generalised length} tells us that in the fusion graph $\Gamma^l$, every vertex $V$ has a path into the component $\Gamma^l_G$. The identity $V\otimes \Ind_G(W)=\Ind_G(\Res_G(V)\otimes W)$ implies that $\Gamma^l_G$ as a graph is isomorphic to the fusion graph of $\Res_G(V)$ as a $kG$-module. If $G$ is finite, $\Res_G(V)$ being faithful implies, by the Burnside-Brauer-Steinberg Theorem, that any simple $kG$-module is reachable from the trivial $kG$-module in  $\Gamma^l_G$. In fact, $\Gamma^l_G$ is strongly-connected, as if $X$ is any $kG$-module, and $X^*$ is its dual, then $X\otimes X^*$ contains a copy of the trivial $kG$-module, so there is a path from any vertex in $\Gamma^l_G$ to the trivial module, and vice versa. 
The same argument in the proof of \cite[Proposition 4.23]{lacabanne2024asymptoticsinfinitemonoidalcategories} now shows that $\Gamma_G^l$ is an FBC which is final in $\Gamma$, so we deduce the statements of \autoref{finite theorem} for $l(n)$ as in the proof of \cite[Theorem 1]{he2025growthproblemsrepresentationsfinite2}. If the group-injective condition is satisfied, or some $V^{\otimes N}$ {has a summand on which $M\setminus G$ acts as zero}, the proof of \cite[Theorem 1]{he2025growthproblemsrepresentationsfinite2} implies that for some large $N'$, $V^{\otimes N'}$ contains a direct summand which is of the form $\Ind_G(P)$, for some projective $kG$-module $P$. This implies there is a path from any vertex in $\Gamma$ to $\Gamma^b_G$. The same argument as above now shows that $\Gamma^b_G$ is an FBC, so the claim follows again as in \cite[Theorem 1]{he2025growthproblemsrepresentationsfinite2}. 

(b). If all simple $kG$-modules are one-dimensional, for the purpose of counting $l(n)$ (which is equal to the sum of all weights of paths of length $n$) we can replace $\Gamma^l_G$ by a single vertex $T$ with self-loop of weight $\dim V$, 
\begin{gather*}
\begin{tikzpicture}[>=stealth, thick]
\node[circle, draw, minimum size=16pt] (v) {T};
\path[->, draw=purple, line width=1.2pt]
(v) edge[loop above, out=135, in=45, min distance=18mm]
node[above, yshift=2pt, text=black] {$\dim V$} (v);
\end{tikzpicture}
\end{gather*}
so that edges ending in any node in $\Gamma^G$ are now replaced by one ending in $T$, with the same label. By the same argument as in (a) above, this is an FBC so the statements of \autoref{finite theorem} for $l(n)$ hold. Moreover, we may assume that $G$ is trivial, so \autoref{finite theorem} gives $l(n)\sim (\dim V)^n$.

(c). By \autoref{all comp factors group}, $V^{\otimes N}$ contains a $kM$-module summand $W$ all of whose composition factors are $kG$-modules. This implies that any $m\in M\setminus G$ acts nilpotently on $W$, and if $m$ is moreover idempotent then it acts as $0$ on $W$. If $M$ is right idempotent generated, then $m=ne$ for some $e\neq 1$ idempotent, so $mW=0$. This is true for all $m\in M\setminus G$, so $W$ is {the desired module on which $M\setminus G$ acts as 0}.

(d). This follows from (a) and (c).
\end{proof}

\begin{Example} We give a few examples.
\begin{enumerate}
\item
Let $FM(T)$ be the free monoid on a finite set $T$. If $V$ is a $kFM(T)$-module with each generator $t\in T$ acting as some singular matrix $N_t$, then (a) of \autoref{inf theorem} gives $l(n)\sim (\dim V)^n$. 
\item Let $M$ be the Temperley--Lieb monoid $TL_m$ or the Motzkin monoid $Mo_m$, then these are regular, and for $V$ on which nonunits act noninvertibly (d) of \autoref{inf theorem} gives $b(n)\sim (\dim V)^n$, for $k$ of any characteristic. This generalizes \cite[Proposition 3C.1]{he2025tensorpowersrepresentationsdiagram} which requires $k=\C$. 
\end{enumerate}
More generally, (d) of \autoref{inf theorem} applies to many important monoids, such as the full transformation monoid, the matrix monoids, and the other diagram monoids in \autoref{fig: table of algebras}. 
\end{Example}

\subsubsection{Rational $kM$-modules}
Rather than considering growth problems for the category of $kM$-modules, we can consider instead the category of rational representations of some linear algebraic monoid. If $M$ is an {irreducible} closed submonoid of the algebraic monoid $M_m$ of $m\times m$ matrices, then it is the Zariski closure of its group of units $G$ in $M_m$. The inclusion $k[M]\xhookrightarrow{} k[G]$ identifies $k[M]$ with the ring of polynomial functions on $G$, and so rational representations of $M$ are the same as polynomial representations of $G$ (see \cite[\S 2.1]{DOTY1998165}): indeed, both are equivalent to comodules over $k[M]$, and this is an equivalence of abelian tensor categories. Thus the growth problem for $M$ {reduces to the growth problem} for $G$, i.e. the case of groups.

\subsubsection{Semigroup case}
The results of the previous section all have analogs where we consider representations of semigroups instead of monoids. If $S$ is a semigroup, a $kS$-module is understood to be non-unital, so that the one-dimensional null representation $Z$ on which $S$ acts as 0 is accepted. The following result in particular proves \cite[Conjecture 4.2]{he2025tensorpowersrepresentationsdiagram} when $S$ is right idempotent generated.

\begin{Theorem}\label{semigroup analogue}
Let $S$ be a semigroup without unit, and suppose that it is dimension-bounded and right-finite with finite set $X$. Let $V$ be a $kS$-module on which no element in $S$ acts invertibly. As above, we write $Z$ for the null representation.
\begin{enumerate}
\item We have that $V^{\otimes \lvert X\rvert }$ has a nonzero submodule on which $S$ acts as 0.
\item We have $l(n)\sim (\dim V)^n$.
\item There is some $N$ such that $V^{\otimes N}$ contains a summand all of whose composition factors are isomorphic to $Z$.
\item If $Z$ is injective, or more generally $V^{\otimes N}$ contains $Z$ as a direct summand for some $N$, then $b(n)\sim (\dim V)^n$.
\end{enumerate}
Suppose moreover that $S$ is right idempotent generated.
\begin{enumerate}[resume]
\item There is some $N$ such that $V^{\otimes N}$ contains $Z$ as a direct summand.
\item We have $b(n)\sim (\dim V)^n$.
\end{enumerate}
\end{Theorem}
\begin{proof}
Similar to the case for monoids. 
\end{proof}

\subsection{Growth problems for virtually abelian groups}

Results from the previous section suggest that the growth problems for a $kM$-module $V$ are often controlled by the $kG$-module $\Res_G(V)$. This motivates us to study the growth problems of virtually abelian groups (or more specifically those with a finite index abelian normal subgroup), which include the groups of units $\widetilde{S}_m$ and $\Z \wr S_m$ of our affine diagram monoids. Growth problems for infinite groups have been studied in \cite{lacabanne2024asymptoticsinfinitemonoidalcategories,coulembier2024fractalbehaviortensorpowers,larsen2024boundsmathrmsl2indecomposablestensorpowers} and are difficult in general. 

\begin{Proposition}\
\begin{enumerate}
\item If $G$ has a subgroup $H$ of finite index $m$ such that all simple modules of $kH$ are bounded in dimension by $r$, then all simple $kG$-modules have dimension $\le mr$ and so \[l(n)\ge \frac{1}{mr}(\dim V)^n.\]
\item In particular, if $G$ has an abelian subgroup of index {$m$} then \[l(n)\ge \frac{1}{m} (\dim V)^n.\]
\end{enumerate}
\end{Proposition}
\begin{proof}
If $V$ is any simple $kG$-module, let $U$ be a simple $kH$-submodule of $V|_H$, then $0\neq \Hom_{kH}(U,V|_{H})\cong \Hom_{kG}(\Ind_H^G(U),V)$, so $\dim V \le \dim \Ind_H^G(U)=mr$.      
\end{proof}

Suppose $G$ has a finite-index normal subgroup $N$. If $V$ is a simple $kG$-module, then the restriction to $N$ of $V$ splits as a direct sum of $g$-conjugates of some simple $kN$-module $U$,

\begin{equation}\label{eqn:clifford}
{V|_N=\bigoplus_{i=1}^t {(U^{g_i})}^{\oplus e}},    
\end{equation}
where if $I_G(U)$ is the inertial subgroup of $U$ (cf. \autoref{semidrect backgroun}) then $t=[ G: I_G(U)]$, {$g_1,\dots, g_t$ is a system of representatives for $G/I_G(U)$} and $e=\dim \Hom (V|_N, U)$. One can therefore say that each $U^{g_i}$ contributes to $1/et$-th of a composition factor of $V$. Note that if $I_G(U)=N$, then {$\Ind_N^G(U)=V$}. (See \cite{Clifford} and also \cite[Theorem 7.1.2]{craven2019representation}.)

\begin{Theorem}\label{virtually abelian thm}
Let $G$ be a group with an abelian normal subgroup $A$ of finite index $m$. Let $V$ be a $kG$-module such that for each $1\neq g\in G/A$, $V|_A$ contains some composition factor $\chi$ such that $g\cdot \chi/\chi$ has infinite multiplicative order. (Where as in \autoref{semidrect backgroun}, {$g\cdot (\chi)(a)=\chi(g^{-1}ag)$.)} Then 
\[l(n)\sim k(n):=\frac{1}{m} (\dim V)^n.\]
\end{Theorem}

\begin{Remark}
The condition is satisfied if $A$ is free abelian, some $\chi$ is faithful and there is $a\in A$ with $C_G(a)=A$. To see this, note that {$(g\cdot \chi/\chi) (a)=\chi(g^{-1}ag a^{-1})$}. Pick $a\in A$ such that {$g^{-1}ag^{} \neq a$}, then {$\chi(g^{-1}aga^{-1})\neq \chi(1)$} by faithfulness. {Thus, for $1\neq g\in G/A$,} $g\cdot \chi/\chi$ must have infinite order: {any non-identity element in a free abelian group has infinite order and $\chi$ is faithful, so $g\cdot\chi/\chi(a)$ (and hence $g\cdot\chi/\chi$) has infinite order.} Note that also this condition forces $G$ to be infinite. 
\end{Remark}
\begin{proof}
  We {work in} the Grothendieck ring where we have $[V^{\otimes n}|_A]=\sum_{w\in J^n} [\chi_w]$, where $J$ is a multiset and $[V|_A]=\sum_{j\in J}[\chi_j]$ is the decomposition into simple $kA$-composition factors (which are all one-dimensional because $k$ is algebraically closed), and for $w=(j_1,\dots,j_n)\in J^n$, $\chi_w=\chi_{j_1}\dots\chi_{j_n}$.

It will suffice to show that as $n\to \infty$, the characters with stabilizer $A$ dominate: by the previous discussion, in this case {each orbit of $m$ such characters corresponds (via induction and restriction) to a single $kG$-composition factor}. {Note that we do not need to assume that $V$ is simple: we are working in the Grothendieck ring so we may replace $[V]$ by the sum of the classes of its composition factors; we then apply Clifford theory to each simple summand.} Let $g_1,\dots, g_k$ be nontrivial coset representatives of $A$ in $G$. For $g_i$, consider a random walk {$R_n=X_1\dots X_n$ (written multiplicatively)} starting at $1$, where the steps {$X_l, l\ge 1$ are the independent and identically distributed random variables} corresponding to multiplication by $g_i\cdot \chi_j/\chi_j$, for $j\in J$ chosen with uniform probability (we can think of $\chi_j$ as some tuple over $k^\times$). Then a fixed point $w\in J^n, g_i\cdot \chi_w=\chi_w$ corresponds to the random walk arriving back at $1$ at step $n$. {Thus it remains to show that the return probability approaches 0 as $n\to \infty$.}

{Fix some $g_i=:g$ and write $\eta_j:=g\cdot \chi_j/\chi_j$, then $H=\langle \eta_j\mid j\in J\rangle$ is a finitely generated abelian group. By assumption, some $\eta\in H$ has infinite order, so we may quotient out by the torsion part of $H$ to obtain a free abelian group $\overline{H}\cong\Z^d, d\ge 1$, in which the image $\overline{\eta}$ of $\eta$ has infinite order. We may project further onto a copy of $\Z$ where the image of $\overline{\eta}$ has infinite order. Let $\pi$ denote the composite of the two projections and write $S_n:=\pi(R_n)=Y_1+\dots+Y_n$ (written additively) then it is enough to show that $\mathbb{P}(S_n=0)\to 0$ as $n\to \infty$. Recall that the \textit{concentration function} of a real random variable $X$ is $\mathcal{Q}_h(X):=\sup_{x\in \R}\mathbb{P}(X\in (x,x+h])$, and for $S_n=Y_1+\dots+ Y_n$ a sum of independent random variables, the Kolmogorov--Rogozin inequality (\cite[Theorem 1]{rogozin1961estimate}) states that for $h>0$, there is $C>0$ such that
\[\mathcal{Q}_h(S_n)\le C\left(\sum_{l=1}^n(1-\mathcal{Q}_h(Y_l))\right)^{-1/2}.\]
Under the present hypotheses, taking $S_n=\pi(R_n)$, we have $\mathcal{Q}_{1/2}(Y_l)=\max_{a\in \Z}\mathbb{P}(Y_l=a)<1$ ($Y_l$ is nonconstant: the case of $Y_l\equiv \pi(\eta)$ is ruled out by observing that $\prod_{j\in J}\eta_j=\prod_j (g\cdot \chi_j)/\chi_j=1$ since $g$ acts on $J$ by permutation, and hence $\mathbb{E}(Y_l)=0$ when switching to additive notation). This implies
$\mathbb{P}(S_n=0)\le \mathcal{Q}_{1/2}(S_n)\le C/\sqrt{n\alpha}$ for some positive $\alpha$. This converges to 0 as $n$ approaches infinity.}

Any element $h\in G\setminus A$ can be written as $g_i a$, for some $i$ and $a\in A$, so if $\chi_w$ is fixed by $h$ then it is fixed by some $g_i$; but the above shows that the proportion of fixed points of each $g_i$ goes to 0. {In other words, in the limit we only need to consider contributions coming from characters $\chi$ with stabilizer $A$. Inducing such a $\chi$ to $G$ results in $m$ conjugates (the maximum possible), see the discussion before the theorem, which implies the result.}
\end{proof}

\begin{Remark}
In \autoref{virtually abelian thm} there is no restriction on the characteristic of $k$, but the condition that $g\cdot \chi_j/\chi_j$ has infinite multiplicative order for some $j$ forces $k$ to have an element of infinite order if $A$ is finitely generated.
\end{Remark}

\begin{Proposition}\label{wreath cor}
If $G=\widetilde{S}_m$ or $\Z\wr S_m$, and $V$ satisfies the hypothesis of \autoref{virtually abelian thm}, then \[l(n)\sim k(n):= \frac{1}{m!}(\dim V)^n.\] 
\end{Proposition}

\begin{proof}
{This follows from \autoref{virtually abelian thm}.}
\end{proof}

\begin{figure}[H]
\centering
\includegraphics[width=0.5\linewidth]{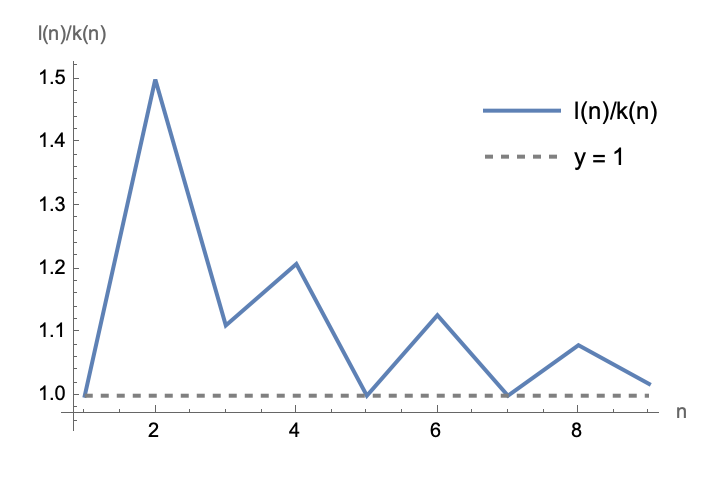}
\caption{The induced representation $\Ind(\chi_{1,2,3})$, a six-dimensional representation of $\Z \wr S_3$. {It is} approximated using a large finite field ($\mathbb{F}_{19699})$, note $19699>3^9$. {The} asymptotic growth rate is $6^{n-1}$.}
\label{fig:ZwrS3}
\end{figure}

We illustrate in \autoref{fig:ZwrS3} the convergence $l(n)/k(n)\to 1$ in the case where $G=\Z\wr S_3$, and $\chi$ is the induction to $\Z \wr S_3$ of the character of $\Z^3$ which takes values $1,2,3$ on the basis vectors of $\Z^3$.

\begin{Example}\label{wallpaper}
Consider the wallpaper group $p_3=\langle x,y,t: xy=yx, t^3=1, txt^{-1}=y, tyt^{-1}=x^{-1}y^{-1} \rangle$, that is $\Z^2 \rtimes C_3$ with $t$ acting by the matrix $\begin{psmallmatrix}
0 & -1 \\
1 & -1.
\end{psmallmatrix}$ This is an index 2 subgroup of $\widetilde{S}_3$. Fix a character $\chi_{a,b}$ of $\Z^2=\langle x,y\rangle$ defined by {$\chi_{a,b}(x)=a, \chi_{a,b}(y)=b$, where $a,b\in k^\times$}, then its conjugates are $\chi_{a,b}^t=\chi_{b,a^{-1}b^{-1}}$ and $\chi_{a,b}^{t^2}=\chi_{a^{-1}b^{-1},a}$. {Let $V=\Ind_{\Z^2}^{p_3}\chi_{a,b}$.} We get \[
V^{\otimes n}\big|_{A}
\;\cong\;
\bigoplus_{k_0+k_1+k_2=n}
\binom{n}{k_0,k_1,k_2}\;
\chi_{\left(
a^{\,2k_0+k_1-n}\,b^{\,k_0+2k_1-n},
\;\;
a^{n-k_0-2k_1}\,b^{\,k_0-k_1}
\right)}.
\]

A character $\chi_{u,v}$ either has stabilizer $\Z^2$ or is fixed also by $C_3$. {Assume that $a,b$ are multiplicatively independent, then the formula above shows that the only $C_3$-fixed character is the trivial one, which occurs only for} $k_0=k_1=k_2=n/3$, and so {\[l(n)=\begin{cases}
3^{n-1} + \frac{2}{3}\binom{n}{n/3, n/3, n/3} &  3\mid n \\
3^{n-1} & \text{else}.
\end{cases}\]} We illustrate the convergence $l(n)/k(n)\to 1$ in \autoref{fig:wallpaper3} below. 
\end{Example}

\begin{figure}[H]
\centering
\includegraphics[width=0.5\linewidth]{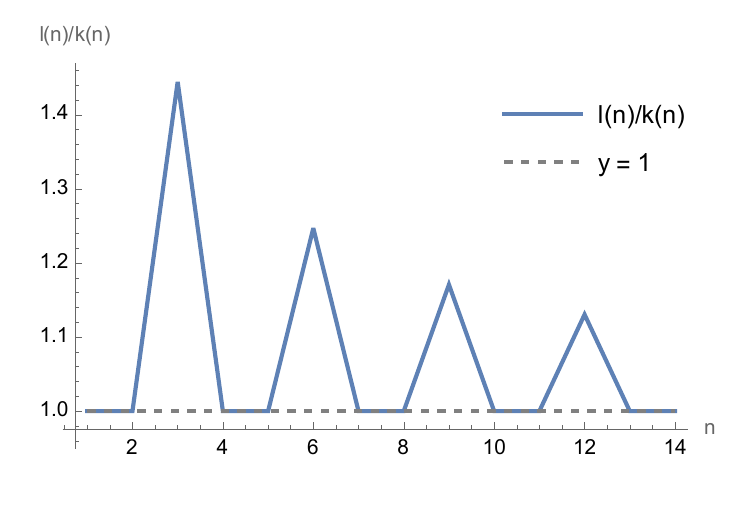}
\caption{An illustration of the convergence $l(n)/k(n)\to 1$, where we take $\chi_{a,b}$ with $a=2, b=3$ in \autoref{wallpaper}.}
\label{fig:wallpaper3}
\end{figure}

\subsection{Growth problems for affine diagram monoids}
We now turn to study the growth problems for the diagram monoids defined in \autoref{alg and monoid section}. 
\begin{Proposition}\label{fin gen nonunit}
Any (reduced) affine or (reduced) periodic monoid associated to a planar diagram monoid or the rook monoid is right-finite.
\end{Proposition} 
\begin{proof}
This follows from the discussions in \autoref{pres_alg}. For example, for the affine or periodic Temperley--Lieb monoid we may take $X=\{e_0,\dots, e_{m-1}\}$. (For planar rook and Motzkin we do not have a presentation, but \autoref{planar rook gen} and \autoref{motzkin gen} give suitable generating sets and relations to move the nonunits to the right.) The claim for the reduced versions follows immediately, since they are quotients of the non-reduced monoids. The affine and periodic rook cases are also clear. 
\end{proof}

We immediately deduce the following:

\begin{Corollary} Let $M$ be any (reduced) affine or (reduced) periodic planar monoid {with group of units $G$}. Let $k$ be any algebraically closed field. If $V$ is a $kM$-module on which no nonunits act invertibly, and such that $\Res_G(V)$ is faithful, then $l(n)\sim (\dim V)^n$. 
\end{Corollary}

\begin{proof}
This follows from (b) of \autoref{inf theorem}, since {the group of units} is either $\Z$ or trivial in this case. 
\end{proof}

\begin{Proposition}\label{sandwich eventually}
Let $M$ be any (reduced) affine or (reduced) periodic monoid associated to a planar diagram monoid or the rook monoid. If $V$ is a $kM$-module on which no nonunits act invertibly, and such that $\Res_G(V)$ is faithful, then there is $C\in \R_{> 0}$ such that $C(\dim V)^n \le l(n)\le (\dim V)^n$, and moreover $r_n(V)\to 1$ as $n\to \infty$ where $r_n(V)$ is as in \autoref{eventually all group}.     
\end{Proposition}
\begin{proof}
It is clear from \autoref{irrep classification theorem} that all simple modules of $M$ have their dimensions bounded above {by} $D$. The results now follow from \autoref{eventually all group}.
\end{proof}



\begin{Proposition}
Let $M$ be a (reduced) periodic monoid associated to a planar monoid. If $V$ is such that no nonunits act invertibly, then $b(n)\sim (\dim V)^n$.     
\end{Proposition}
\begin{proof}
The monoids are all right-finite, dimension-bounded, and have trivial group of units, so by (d) of \autoref{inf theorem} it suffices to show that all of these monoids are right idempotent generated. They are right idempotent generated for the same (diagrammatic) reason that the ordinary planar diagram monoids are ({see \cite[Lemma 4C.1]{KhSiTu-monoidal-cryptography}): namely, for a diagram $D$,  as before let $D^\star$ denote its reflection across a horizontal axis, then $DD^\star D=D$ after removing closed components, see the diagram below.} 
\[\begin{tikzpicture}

\begin{scope}[x=0.55cm,y=0.55cm]
  \draw[usual] (2,3) arc[start angle=180,end angle=360,x radius=.5,y radius=.3];

  \draw[usual] (0,3)--(2,2);
  \draw[usual] (1,3)--(3,2);

  \draw[usual] (0,2) arc[start angle=180,end angle=0,x radius=.5,y radius=.3];
  \draw[usual] (0,2) arc[start angle=180,end angle=360,x radius=.5,y radius=.3];

  \draw[usual] (2,2)--(0,1);
  \draw[usual] (3,2)--(1,1);

  \draw[usual] (2,1) arc[start angle=180,end angle=0,x radius=.5,y radius=.3];
  \draw[usual] (2,1) arc[start angle=180,end angle=360,x radius=.5,y radius=.3];

  \draw[usual] (0,1)--(2,0);
  \draw[usual] (1,1)--(3,0);

  \draw[usual] (0,0) arc[start angle=180,end angle=0,x radius=.5,y radius=.3];
\end{scope}

\node at (2.45,1.0) {$=$};


\begin{scope}[shift={(3.15,0.5)},x=0.55cm,y=0.55cm]
  \draw[usual] (2,1) arc[start angle=180,end angle=360,x radius=.5,y radius=.3];
  \draw[usual] (0,1)--(2,0);
  \draw[usual] (1,1)--(3,0);
  \draw[usual] (0,0) arc[start angle=180,end angle=0,x radius=.5,y radius=.3];
\end{scope}

\node at (0.82,-0.42) {$DD^{*}D$};
\node at (3.95,-0.42) {$D$};

\end{tikzpicture}\]
{Now if we let $e=D^\star D$, then this is an idempotent and $D=De$. It is clear that $e$ is not identity if $D$ is not. This shows $M$ is right idempotent generated.}
\end{proof}

\begin{figure}[H]
\centering
\includegraphics[width=0.5\linewidth]{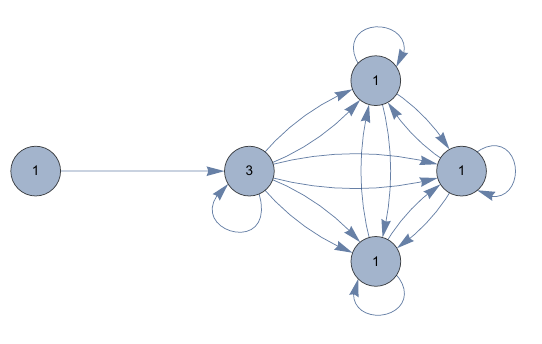}
\caption{The fusion graph of $\Delta(1,1)$. The vertices are labelled with their dimensions.}
\label{fig:W31}
\end{figure}

We illustrate in \autoref{fig:W31} the fusion graph $\Gamma^b$ of the cell module $\Delta(1,1)$ corresponding to $(\lambda,z)=(1,1)$ for $\overline{aTL}_3$. Note that $\Delta(1,1)$ is in fact a module for the monoid underlying an `uncoiled' algebra (see \cite[\S 3.1]{langlois2023uncoiled}), whose group of units is the cyclic group of order 3. This corresponds to the three one-dimensional $kG$-modules forming a strongly-connected component $\Gamma^b_G$ in \autoref{fig:W31}.

\subsubsection{Growth problems for symmetric affine diagram monoids}

Recall from \autoref{fin dim section} that many of the simple modules for the reduced affine symmetric diagram algebras may be obtained via inflation from their $r$-reduced quotients, which are finite-dimensional. {Throughout this section we assume $k=\C$} and study $b(n)$ for these modules. 

The group of units of the monoids $aPa_{m,r},aRoBr_{m,r}$ and $aBr_{m,r}$ are the generalized symmetric groups $C_r\wr S_m = C_r^m \rtimes S_m$ (= complex reflection groups $G(r,1,m)$), which we denote by $S(r,m)$. {As in} \autoref{rep of wreath}, the simple $\C S(r,m)$-modules are parameterized by $r$-multipartitions ${\lambda}=(\lambda^{(1)},\dots, \lambda^{(r)})$ with $\lvert \lambda^{(1)}\rvert+\dots + \lvert \lambda^{(r)}\rvert = m$. In more detail, write $r_i=\lvert \lambda^{(i)}\rvert\ge 0$ and let $S^{\lambda^{i}}$ be the corresponding Specht module. Then all simple $\C S(r,m)$-modules can be constructed as:
$$W_{\lambda}=\Ind_K^{S(r,m)} \boxtimes_{i=1}^r (\chi_i^{\boxtimes r_i}\boxtimes S^{\lambda^{i}}), $$
where $\boxtimes$ denotes the {external tensor} product, $K$ is the \textit{inertia group}, $K:=\prod_{i=1}^r (C_r^{r_i}\rtimes S_{r_i})$, and $\chi_i$ is a character of $C_r$. In particular, we get 
$$\dim W_{{\lambda}}=\frac{m!}{\prod_{i=1}^r r_i!} \prod_{i=1}^r f^{\lambda^{i}},$$
where $f^{\lambda^{i}}= \dim S^{\lambda^{i}}$.

We first study the growth problems for $S(r,m)$.
\begin{Proposition}\label{generalised sym}
If $V$ is a faithful simple representation of $S(r,m)$, then 
\[l(n)=b(n)\sim a(n)=\begin{cases}
(\sum_{k=0}^{\floor{m/2}}\frac{1}{r^k (m-2k)!k!2^k}+\frac{1}{(m/2)! r^m}(\frac{r}{2})^{m/2}\cdot (-1)^n)\cdot (\dim V)^n & $r,m$ \ \text{even},\\
\sum_{k=0}^{\floor{m/2}}\frac{1}{(m-2k)!k!2^k r^k}\cdot (\dim V)^n & \text{else}.
\end{cases}\]
\end{Proposition}
\begin{proof}
We apply \autoref{finite theorem} (to the case $M=G$) and note that $Z_V(G)=Z(G)$ because $V$ is simple and faithful. 
Write $g$ for the generator of $C_r$, then the center of $S(r,m)$ is $\{z_t:=(g^t,\dots, g^t; 1): 0\le t\le r-1\}\cong C_r$. {Write $\chi_\lambda$ for the character of $W_\lambda$, as defined above.} The sum over the column in the character table of $S(r,m)$ corresponding to $z_t$ is 
\[\begin{aligned}
S_t
&= \sum_{\lambda}\chi_{\lambda}(z_t) = \sum_{r_1+\dots+r_r=m}
\zeta^{t\sum_{i=1}^r (i-1)r_i}\,
\frac{m!}{\prod_{i=1}^r r_i!}
\sum_{\lambda^{(1)}\vdash r_1}\cdots
\sum_{\lambda^{(r)}\vdash r_r}
\prod_{i=1}^r f^{\lambda^{(i)}} \\
&= \sum_{r_1+\dots+r_r=m}
\zeta^{t\sum_{i=1}^r (i-1)r_i}\,
\frac{m!}{\prod_{i=1}^r r_i!}
\prod_{i=1}^r\left(\sum_{\lambda^{(i)}\vdash r_i} f^{\lambda^{(i)}}\right) \\
&= \sum_{r_1+\dots+r_r=m}
\zeta^{t\sum_{i=1}^r (i-1)r_i}\,
m!\prod_{i=1}^r \frac{S(r_i)}{r_i!},
\end{aligned}\]
where $\zeta$ denotes a primitive $r$-th root of unity, $\zeta^{t\sum_{i=1}^r(i-1)r_i}\cdot \dim W_{\lambda}=\chi_{\lambda}(z_t)$, and $S(k)$ is the sum of the dimensions of all simple $\C S_k$-modules.
Let $F(x):=\sum_{k\ge 0} {S(k)}x^k/{k!} = e^{x+x^2/2} $
be the exponential generating function for the number of involutions in $S_k$, see \cite[A000085]{OEIS}. Then we have \[S_t=m! [x^m]\prod_{i=0}^{r-1} F(\zeta^{it}x)= m! [x^m]\exp\Big( \sum_{i=0}^{r-1} \bigl( \zeta^{ti}x + (\zeta^{ti}x)^2/2 \bigr) \Big).
\]
From the above equation we deduce that \[S_0=
\sum_{k=0}^{\lfloor m/2 \rfloor} \frac{m!}{(m-2k)!k!2^k} r^{m-k}, \quad S_t=\frac{m!}{(m/2)!}\big(r/2 \big)^{m/2}\] if {$r$ and $m$} are both even and $t=r/2$ (otherwise, $S_t=0$). The formula for $a(n)$ follows from \autoref{finite theorem}, noting that $\lvert S(r,m)\rvert = r^m\cdot m!$. Finally, note that $l(n)=b(n)$ since $\C G$ is semisimple.     
\end{proof}

\begin{Remark}
\autoref{generalised sym} generalizes the result of \cite[Example 2.3]{coulembier2023asymptotic} for the symmetric group, and \cite[Example 4]{lacabanne2023asymptotics}.
\end{Remark}

\begin{Remark}
If we let $r\to \infty$ in the formula for $a(n)$ in \autoref{generalised sym}, we get $(1/m!)\cdot (\dim V)^n,$ which is $l(n)$ for $\Z \wr S_m$ (under the right assumptions on $V$) by \autoref{wreath cor}. Of course, for an infinite group over $\C$ we no longer have semisimplicity, so $l(n)\neq b(n)$.
\end{Remark}

\begin{Lemma}\label{group injective lemma}
The monoids $aPa_{m,r},aRoBr_{m,r}$ and $aBr_{m,r}$ satisfy the group-injective condition over $\C$.    
\end{Lemma}
\begin{proof}
By \cite[Proposition 3B.4; Proposition 3B.6]{he2025tensorpowersrepresentationsdiagram}, the monoids $N=Br_m, RoBr_m, Pa_m$ possess idempotents $e$ such that $Ne$ is isomorphic to $\Ind_G(V_{sgn})$, where $V_{sgn}$ is the sign representation of {$S_m$}. Denote the corresponding $r$-reduced affine monoid by $M$. In all cases, $e\in \C N\subseteq \C M$ may be written as a linear combination of only elements of {$S_m$} and the `stops' $p_i$ (cf. \autoref{rook brauer section}). Let $f:=(1/r^m) \cdot \sum_{g \in C_r^m} g$ be {the primitive idempotent for the trivial representation of $C_r^m$ (which we identify with the subgroup generated by the windings)}, then $hfh^{-1}=f$ for all $h\in S_m$, whence $ef=fe$ in $\C M$ since the stops also commute with $C_r^m$. It is clear that $ef$ is an idempotent on which the group of units $C_r^m \rtimes S_m$ acts as scalar, and any nonunit acts as 0. {Moreover, $ef\neq 0$: indeed, in the quotient $\C S(r,m)$ of $\C M$, $ef$ is the primitive idempotent corresponding to the simple representation $W$ which is trivial on $C_r^m$ and is $V_{sgn}$ on $S_m$.} In other words, $\C Mef \cong \Ind_G(W)$ {is projective as a $\C M$-module, for $W$ a projective (= simple) $\C S(r,m)$-module.} Recall that $M$ possesses an involution $\iota$ that flips diagrams across a horizontal axis. If $V$ is a $\C M$-module, we can define $D(V)=\textup{Hom}(V,\C)$ with the action $x\cdot f(v)=f(\iota(x) v)$. We then get a {simple-preserving} duality $V\mapsto D(V)$ sending projective modules to injective modules {(see e.g. \cite[\S 2.3]{belletete2018restriction} for the argument)}. {Note that $W$ is self-dual as a $\C S(r,m)$-module, and hence $\C Mef$ is self-dual with respect to $D$: $ D\bigl(\operatorname{Ind}_G(W)\bigr)
\cong \operatorname{Ind}_G(W^*)
\cong \operatorname{Ind}_G(W)$, where $(-)^*$ denotes the usual duality of group representations.} Thus $\C Mef\cong D(\C Mef)$ is also injective. {This completes the proof.}   
\end{proof}

\begin{Proposition}
Let $M= aPa_{m,r},aRoBr_{m,r}$ or $aBr_{m,r}$. If $V$ is a $\C M$-module on which no nonunit acts invertibly {and whose restriction to the group of units is faithful and simple}, then 
\[l(n)\sim b(n)\sim \begin{cases}
(\sum_{k=0}^{\floor{m/2}}\frac{1}{r^k (m-2k)!k!2^k}+\frac{1}{(m/2)! r^m}(\frac{r}{2})^{m/2}\cdot (-1)^n)\cdot (\dim V)^n & $r,m$ \ \text{even},\\
\sum_{k=0}^{\floor{m/2}}\frac{1}{(m-2k)!k!2^k r^k}\cdot (\dim V)^n & \text{else}.
\end{cases}\]
\end{Proposition}
\begin{proof}
Immediate from \autoref{group injective lemma} and \autoref{finite theorem}.
\end{proof}

\bibliographystyle{alphaurl}
\bibliography{sources}

\end{document}